\documentclass[aop,seceqn,MSNbibl,citesort,dvips]{arximspdf}
\usepackage{graphics}
\doi{10.1214/07-AOP388}
\volume{36}
\issue{6}
\pubyear{2008}
\firstpage{2235}
\lastpage{2279}

\makeatletter

\newproclaim{remark}{Remark}[section]
\newtheorem{proposition}[remark]{Proposition}
\newtheorem{theorem}[remark]{Theorem}

\newproclaim{step}{Step}
\newproclaim{stepp}{Step}

\newtheorem{theoremm}{Theorem}[section]
\newtheorem{lemmaa}[theoremm]{Lemma}
\newproclaim{remarkk}[theoremm]{Remark}

\makeatother

\begin{document}
\begin{frontmatter}

\title{A tree approach to $p$-variation and to integration}
\runtitle{A tree approach to $p$-variation and to integration}
\begin{aug}
\author[A]{\fnms{Jean} \snm{Picard}\corref{}\ead[label=e1]{Jean.Picard@math.univ-bpclermont.fr}}
\runauthor{J. Picard}
\affiliation{Universit\'{e} Blaise Pascal}
\address[A]{Laboratoire de Math\'{e}matiques (CNRS UMR 6620)\\
Universit\'{e} Blaise Pascal\\
63177 Aubi\`{e}re Cedex\\
France\\
\printead{e1}} %adresu isvedimo komanda gale!
\end{aug}

% HISTORY:
\received{\smonth{5} \syear{2007}}
\revised{\smonth{12} \syear{2007}}

% ABSTRACT
%
\begin{abstract}
We consider a real-valued path; it is possible to associate a tree to this
path, and we explore the relations between the tree, the properties of
\mbox{$p$-variation} of the path, and integration with respect to the path. In
particular, the fractal dimension of the tree is estimated from the
variations of the path, and Young integrals with respect to the path, as
well as integrals from the rough paths theory, are written as integrals on
the tree. Examples include some stochastic paths such as martingales, L\'{e}vy
processes and fractional Brownian motions (for which an estimator of the
Hurst parameter is given).\looseness=-1
\end{abstract}

% KEYWORDS
%
\pdfkeywords{60G17, 60H05, 26A42,
Lebesgue--Stieltjes integrals,
rough paths, real trees,
variations of paths,
fractional Brownian motion,
Levy processes}

\end{frontmatter}

%s1 ###
\section{Introduction}\label{sec1}

Consider a continuous path $\omega\dvtx[0,1]\to\mathbb{R}$. The
$p$-variation of
$\omega$ is defined for $p\ge1$ by
\[
V_p(\omega):=\sup_{(t_i)}\sum_i|\omega(t_{i+1})-\omega
(t_i)|^p
\]
for subdivisions $(t_i)$ of $[0,1]$. It is well known that the
finiteness of
$V_p(\omega)$ is closely related to the possibility of constructing
integrals $\int_0^1\rho\,d\omega$ for some functions $\rho$. The simplest
case is when $V_1(\omega)$ is finite ($\omega$ has finite variation);
then a
signed measure $d\omega=d\omega^+-d\omega^-$ (the Lebesgue--Stieltjes
measure) is defined from $\omega$, and the integral is well defined
for any
bounded Borel function $\rho$; if moreover~$\rho$ has left and right limits,
then the integral is also a Riemann--Stieltjes integral (it is the limit of
Riemann sums). If now $\omega$ has infinite variation ($V_1(\omega
)=\infty$)
but $V_p(\omega)$ is finite for a larger value of $p$, it was proved by
Young \cite{young36} that a Riemann--Stieltjes integral can still be
constructed as soon as $V_q(\rho)$ is finite for $q$ such that $1/p+1/q>1$;
as an application, one can consider and solve stochastic differential
equations driven by a multidimensional path with finite $p$-variation if
$p<2$ (in particular a typical fractional Brownian path with Hurst parameter
$H>1/2$). If now $p$ is greater than 2, Lyons's theory of rough paths
\cite{lejay03,lyons98,lyonscl07,lyonsqian02} provides a richer framework
which is still suitable to consider and solve these equations.

On the other hand, one can associate to $\omega$ a metric space
$(\mathbb{T},\delta)$ which is a compact real tree and which can be
used to
describe the excursions of $\omega$ above any level; see
\cite{duquesne07,duqleg05} or Chapter 3 of \cite{evans07}. The tree
$\mathbb{T}$
can be endowed with its length measure $\lambda$, and our aim is to relate
the properties of $(\mathbb{T},\delta,\lambda)$ to the questions of
$p$-variation
of $\omega$ and of integration with respect to $\omega$. These
questions are
also considered for c\`{a}dl\`{a}g paths $\omega$ (paths which are right-continuous
and have left limits), since these paths can be considered as time-changed
continuous paths. As an application, we consider the case where $\omega
$ is
a path of a stochastic process such as a L\'{e}vy process or a fractional
Brownian motion (the case of a standard Brownian path has been
considered in
\cite{picard06}).

In Section \ref{trees}, we introduce the tree $\mathbb{T}$ and study
its basic
properties. In particular, in the finite variation case, we work out the
interpretation of its length measure $\lambda$ by means of the
Lebesgue--Stieltjes measure of $\omega$, extending a result of
\cite{duquesne07}; this result is fundamental for the construction of
integrals in Section \ref{integrals} (see below). We also explain how the
tree can be defined in the c\`{a}dl\`{a}g case.

In Section \ref{pvariation}, we see in Theorem \ref{dimtree} (Theorem
\ref{vpomp} for the c\`{a}dl\`{a}g case) that the finiteness of $V_p(\omega)$ is
related to some metric properties of $\mathbb{T}$, particularly its
upper box
dimension $\operatorname{\overline{dim}}\mathbb{T}$; more precisely,
%
%e1.1 ###
\begin{equation}\label{vpomega}
\cases{
\displaystyle V_p(\omega)=\infty,&\quad if $1\le p<\operatorname{\overline{dim}}\mathbb
{T}$,\cr
\displaystyle V_p(\omega)<\infty,&\quad if $p>\operatorname{\overline{dim}}\mathbb{T}$.}
\end{equation}
We give applications of these results to martingales, fractional Brownian
motions and L\'{e}vy processes. We prove in particular that upper box and
Hausdorff dimensions of $\mathbb{T}$ coincide for fractional Brownian motions
(with Hurst parameter $H$) and stable L\'{e}vy processes (with index
$\alpha$);
we also construct an estimator of $H$ based on $\mathbb{T}$, which can be
computed by means of a sequence of stopping times (Proposition~\ref{estimator}).

The aim of Section \ref{integrals} is to construct integrals with
respect to
$\omega$ by means of the tree. Let us assume that $\omega$ is
continuous and
$\omega(0)=\omega(1)=\inf\omega$ (considering the general case adds some
notational complication). The construction of the integral is based on the
following remark (Propositions \ref{lebesgue} and \ref{finitevar}): when
$\omega$ has finite variation, the positive and negative parts
$d\omega^+$
and $d\omega^-$ of $d\omega$ can be viewed as the images of the length
measure $\lambda$ by two maps $\tau\mapsto\tau^\nearrow$ and
$\tau\mapsto\tau^\nwarrow$ from $\mathbb{T}$ to $[0,1]$; thus
%
%e1.2 ###
\begin{equation}\label{annonce}
\int_0^1\rho\,d\omega=\int_\mathbb{T}\bigl(\rho(\tau^\nearrow
)-\rho
(\tau
^\nwarrow)\bigr)
\lambda(d\tau).
\end{equation}
When $\omega$ has infinite variation, this procedure can still be
applied to
construct $d\omega^+$ and $d\omega^-$; these measures are $\sigma$-finite
but no more finite. However, \textup{(\ref{annonce})} can be viewed
as a definition
of $\int\rho\,d\omega$ provided the term in the right-hand side is
integrable; this means that the tree can provide a mechanism by means of
which $d\omega^+$ and $d\omega^-$ compensate each other. For
instance, if
$1/p+1/q>1$,
\[
V_p(\omega)<\infty, \qquad V_q(\rho)<\infty
\quad\Longrightarrow\quad\int_\mathbb{T}|\rho(\tau^\nearrow)-\rho
(\tau
^\nwarrow
)|
\lambda(d\tau)<\infty.
\]
Moreover, in this case, the integral defined by \textup{(\ref
{annonce})} coincides
with the Young integral (Theorems \ref{young} and \ref{cadlag}).
Consequently, differential equations driven by multidimensional paths with
finite $p$-variation with $p<2$ enter our framework. Actually, we may take
$p>2$ for one of the components (Theorem \ref{edo}); this is due to
the fact
that the condition $V_q(\rho)<\infty$ can be replaced by some weaker
condition $V_q(\rho|\omega)<\infty$. We also prove that the tree approach
can be used to consider multidimensional fractional Brownian motions with
parameter $H>1/3$ (Theorem \ref{intfbm}); in this case, the right-hand side
of \textup{(\ref{annonce})} should be understood as a generalized
integral on
$\mathbb{T}$
(a limit of integrals on subtrees $\mathbb{T}^a$ obtained by trimming
$\mathbb{T}$),
and we recover the integrals of the rough paths theory.

The \hyperref[mixing]{Appendix} is devoted to two results which are needed in the
article, and which may also be of independent interest. In Appendix
\ref{mixing}, we prove that increments of fractional Brownian motions are
asymptotically independent from the past. In Appendix \ref
{roughpaths}, we
study the time discretization of integrals in the rough paths calculus,
in a
spirit similar to \cite{feyelprad06,gubin04}.

\begin{remark}
A lot of work has been devoted to the links between random trees and
excursions of some stochastic processes; these links are an extension
of the
classical Harris correspondence between random walks and random finite
trees. Historically, they have first been investigated in the context of
Brownian excursions in \cite{nevpit89,legall91,aldous93} (see also the
courses \cite{pitman06,evans07}) with the aim of studying branching
processes. In order to consider more general branching mechanisms, L\'{e}vy
trees, defined by means of L\'{e}vy processes $X$ without negative jumps, have
been introduced and studied in \cite{leglejan98a,duqleg02}; they
have been
related to the notion of real tree in \cite{duqleg05}. However, we
will not
focus here on properties of L\'{e}vy trees; a L\'{e}vy tree is indeed a tree which
is associated to some continuous process related to $X$ (the height
process), whereas we will rather consider in our applications the tree which
is associated directly to the L\'{e}vy process $X$.
\end{remark}

\begin{remark}
We work out here a nonlinear approach to integration with respect to
one-dimensional paths; consequently, the integral with respect to
$\omega_1+\omega_2$ is not simply related to integrals with respect to
$\omega_1$ and $\omega_2$; moreover, integration with respect to a
multidimensional path can be worked out by summing integrals with
respect to
each component, but this depends on the choice of a frame.
\end{remark}

\begin{remark}
In the proofs of this article, the letter $C$ will denote constant numbers
which may change from line to line. For quantities depending on the
path~$\omega$ of a stochastic process, we will rather use the notation
$K=K(\omega)$.
\end{remark}

%s2 ###
\section{Paths and trees}
\label{trees}

In this section, we first define the tree associated to a continuous path,
describe its length measure, and extend these objects to c\`{a}dl\`{a}g paths.

%s2.1 ###
\subsection{Basic definitions and properties}\label{basic}

Consider a continuous function $(\omega(t); 0\le t\le1)$. The function
%
%e2.1 ###
\begin{equation}\label{delta}
\delta(s,t):=\omega(s)+\omega(t)-2\inf_{[s,t]}\omega
\end{equation}
is a semi-distance on $[0,1]$, where
\[
\delta(s,t)=0\quad\Longleftrightarrow\quad\omega(s)=\omega(t)=\inf
_{[s,t]}\omega.
\]
The quotient metric space $\mathbb{T}=([0,1]/\delta,\delta)$ is a real
tree; this
means that between any two points $\tau_1$ and $\tau_2$ in $\mathbb{T}$,
there is
a unique arc denoted by $[\tau_1,\tau_2]$ ($\mathbb{T}$ is a
topological tree),
and that $[\tau_1,\tau_2]$ is isometric to the interval
$[0,\delta(\tau_1,\tau_2)]$ of $\mathbb{R}$; see \cite{duqleg05}.
Actually, real
trees can also be characterized as connected metric spaces satisfying the
so-called four-point condition, and one can use this condition to prove that
$\mathbb{T}$ is a real tree; see \cite{duquesne07,evans07}. We will
denote by
$\pi$ the projection of $[0,1]$ onto $\mathbb{T}$; notice that if
$\omega
$ is
constant on some interval $[s,t]$, then all the points of this interval are
projected on the same point of $\mathbb{T}$. The continuity of $\pi$
follows from
the continuity of $\omega$; in particular, $\mathbb{T}$ is compact.
In this
article we implicitly assume that $\omega$ is not constant, so that
$\mathbb{T}$
is not reduced to a singleton.

%
%f1 ###
\begin{figure}[b]

\includegraphics{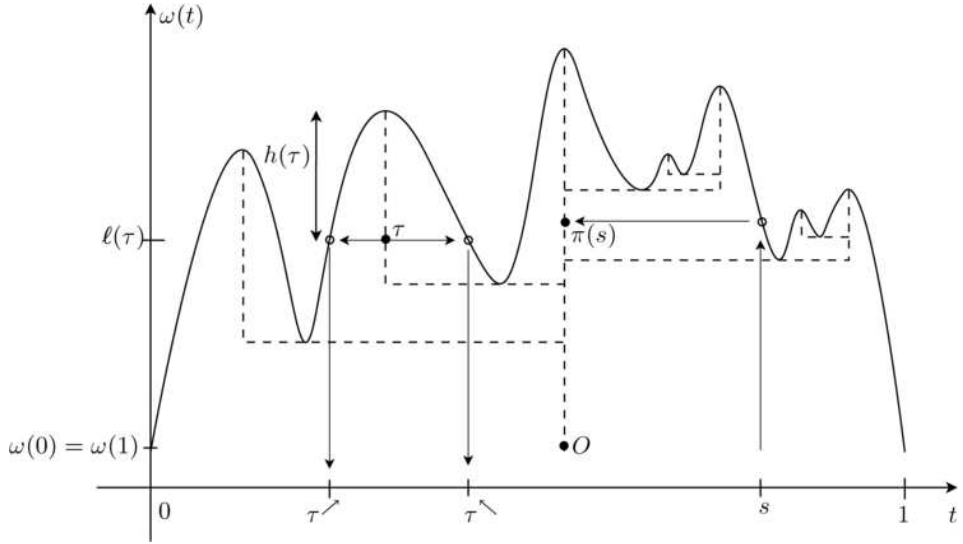}

\caption{An example of path $\omega$ with its tree $\mathbb{T}$
represented by
dashed lines (the vertical lines represent points of the skeleton, and each
branching point is represented by a horizontal line); maps $\tau
\mapsto
\tau^\nearrow$,
$\tau\mapsto\tau^\nwarrow$ and $s\mapsto\pi(s)$ are also
depicted.}\label{fig1}
\end{figure}

We now suppose $\pi(0)=\pi(1)$, or equivalently
%
%e2.2 ###
\begin{equation}\label{om01}
\omega(0)=\omega(1)=\inf_{[0,1]}\omega.
\end{equation}
An example is given in Figure \ref{fig1}. We explain at the end of the subsection how general paths can be
reduced to
this case. Under this condition, $\mathbb{T}$ becomes a rooted tree by
considering \mbox{$\pi(0)=\pi(1)=O$} as the root of the tree, and we can say
that a
point $\tau_1$ is above $\tau_2$ if $\tau_2\in[O,\tau_1]$.

We consider on $\mathbb{T}$ the level function $\ell$ defined by
%
%e2.3 ###
\begin{equation}\label{ell}
\ell(\tau):=\omega(0)+\delta(O,\tau).
\end{equation}
Then $\omega=\ell\circ\pi$. For $\tau$ in $\mathbb{T}$, define
\[
\tau^\nearrow:=\inf\pi^{-1}(\tau), \qquad \tau^\nwarrow:=\sup\pi
^{-1}(\tau),
\]
so that
\[
\omega(\tau^\nearrow)=\omega(\tau^\nwarrow)=\inf_{[\tau
^\nearrow,\tau
^\nwarrow]}\omega
=\ell(\tau).
\]
In particular $O^\nearrow=0$ and $O^\nwarrow=1$. The set
$\pi([\tau^\nearrow,\tau^\nwarrow])$ is exactly the set of points above
$\tau$. If now we consider the set
$\pi([\tau^\nearrow,\tau^\nwarrow])\setminus\{\tau\}$ of points
which are
strictly above $\tau$, it is made of connected components which are
subtrees, and which are called the branches above $\tau$; each of these
branches is the projection of a connected component of
$[\tau^\nearrow,\tau^\nwarrow]\setminus\pi^{-1}(\tau)$, and
corresponds to
an excursion of $\omega$ above level $\ell(\tau)$. If there is more
than one
branch above $\tau$, then $\tau$ is said to be a branching point;
this means
that there is more than one excursion, and the times between these
excursions are local minima of $\omega$ (a local minimum may be a constancy
interval). On the other hand, if there is no branch above $\tau$, then
$\tau$ is said to be a leaf; this means that
$\pi^{-1}(\tau)=[\tau^\nearrow,\tau^\nwarrow]$, so this holds when
$\tau^\nearrow=\tau^\nwarrow$ or when $[\tau^\nearrow,\tau
^\nwarrow]$
is a
constancy interval of $\omega$. Local maxima of $\omega$ are
projected on
leaves of $\mathbb{T}$, but there may be leaves which are not
associated to local
maxima. Points which are not leaves constitute the skeleton $S(\mathbb
{T})$ of
the tree.

We say that $\omega$ is piecewise monotone if there exists a finite
subdivision $(t_i)$ of $[0,1]$ such that $\omega$ is monotone on each
$[t_i,t_{i+1}]$. We also say that $\mathbb{T}$ is finite if it has
finitely many
leaves. If $\mathbb{T}$ is not finite, then it has infinitely many branching
points, or it has at least a branching point with infinitely many branches
above it; in both of these cases, $\omega$ has infinitely many local minima
and is therefore not piecewise monotone. Conversely, if $\omega$ is not
piecewise monotone, then it has infinitely many local maxima, and each of
them is projected on a different leaf of $\mathbb{T}$, so $\mathbb
{T}$ is not finite.
Thus
%
%e2.4 ###
\begin{equation}\label{piecewise}
\mbox{$\omega$ is piecewise monotone}\quad\Longleftrightarrow\quad
\mbox{$\mathbb{T}$ is finite.}
\end{equation}

We shall also need an operation called trimming, or leaf erasure, due to
\cite{neveu86} (see also \cite{kesten86,nevpit89,evanspitwin06,evans07});
to this end, we introduce the function
%
%e2.5 ###
\begin{equation}\label{htau}
h(\tau):=\sup\{\omega(t)-\ell(\tau);\tau^\nearrow\le t\le
\tau^\nwarrow\}.
\end{equation}
This is the height of the (or of the highest) branch above $\tau$. In
particular, $h(\tau)=0$ if and only if $\tau$ is a leaf.

%
%f2 ###
\begin{figure}

\includegraphics{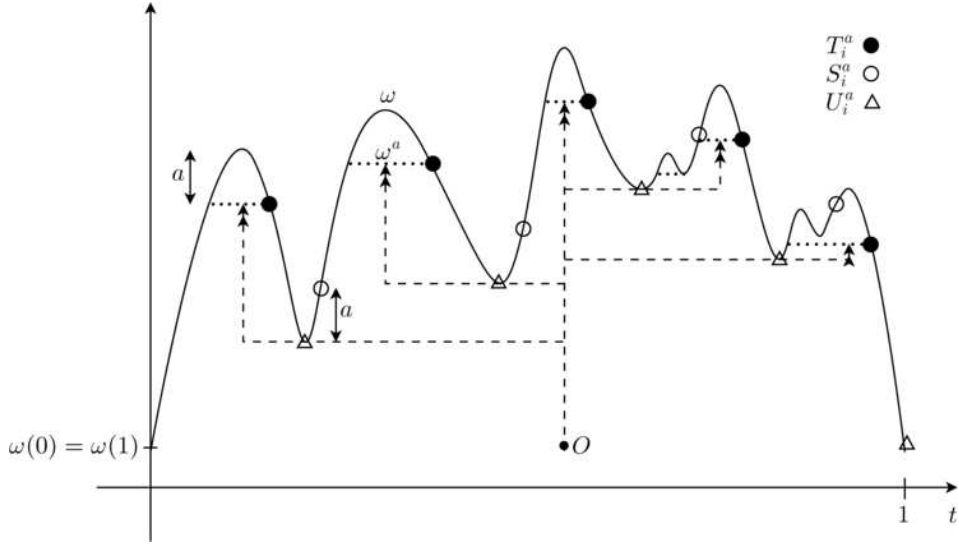}

\caption{The trimmed tree $\mathbb{T}^a$ is represented by dashed
lines, and
its leaves
by double arrows; the flattened path $\omega^a$ is represented by dots
when it differs from
$\omega$; times $T_i^a$, $S_i^a$ and $U_i^a$, $i\ge1$, are respectively
represented
on the curve by bullets, circles and triangles.}\label{fig2}
\end{figure}

Now consider the trimmed tree
%
%e2.6 ###
\begin{equation}\label{deftreea}
\mathbb{T}^a:=\{\tau\in\mathbb{T};h(\tau)\ge a\}.
\end{equation}
Then $\mathbb{T}^a$ is nonempty if and only if
$\|\omega\|:=\sup\omega-\inf\omega\ge a$, and in this case, it is
a rooted
subtree of $\mathbb{T}$ (it contains the root $O$). An example is drawn in Figure \ref{fig2}. As $a\downarrow
0$, the tree
$\mathbb{T}^a$ increases to the skeleton of $\mathbb{T}$; each branch
grows at unit
speed, and a new branch appears at $\tau$ if $\tau$ is a branching
point of
$\mathbb{T}$ such that one of the branches above $\tau$ has height
exactly $a$,
and another one has height at least $a$. This subtree has been
introduced in
\cite{nevpit89} and is related to $a$-minima and $a$-maxima of the path.
More precisely, starting with $S_0^a=T_0^a=0$, define
%
%e2.7 ###
\begin{equation}\label{tiasia}
\cases{\displaystyle T_{i+1}^a:=\inf\biggl\{t\in
[S_i^a,1]; \omega(t)-\sup_{[S_i^a,t]}\omega<-a\biggr\},\cr
\displaystyle S_{i+1}^a:=\inf\biggl\{t\in
[T_{i+1}^a,1]; \omega(t)-\inf_{[T_{i+1}^a,t]}\omega>a\biggr\},\cr
\displaystyle N^a:=\inf\{i; T_i^a\mbox{ or }S_i^a=\inf\varnothing\}.}
\end{equation}
Actually, in the case $\pi(0)=\pi(1)$, $T_{N^a}^a$ is still well defined,
but not $S_{N^a}^a$ (notice in particular that if $\omega$ is a path
of an
adapted stochastic process, then $S_i^a$ and $T_i^a$ are stopping times).
Then $N^a$ is the number of leaves of $\mathbb{T}^a$; the set of leaves
$\partial\mathbb{T}^a$ and the set of times $(T_i^a; 1\le i\le N^a)$
are in
bijection by means of $\pi$ and its inverse map $\tau\mapsto\tau
^\nwarrow$.
Moreover
%
%e2.8 ###
\begin{equation}\label{inftiti}
\inf_{[T_i^a,T_{i+1}^a]}\omega=\inf_{[T_i^a,S_i^a]}\omega
=\omega(S_i^a)-a \qquad\mbox{for $1\le i<N^a$.}
\end{equation}

The approximation of $\mathbb{T}$ by $\mathbb{T}^a$ can also be
interpreted as an
approximation of the path $\omega$; trimming the tree is equivalent to
flattening some excursions of the path. More precisely, let $\pi^a(t)$ be
the projection of $\pi(t)$ on $\mathbb{T}^a$ (assuming $\mathbb
{T}^a\ne
\varnothing$),
and let
%
%e2.9 ###
\begin{equation}\label{defomegaa}
\omega^a=\ell\circ\pi^a
\end{equation}
for the level function $\ell$ defined in \textup{(\ref{ell})}. Then
$\mathbb{T}^a$
is the
associated tree of $\omega^a$. The path $\omega^a$ is continuous, is
obtained from $\omega$ by means of the change of time
\[
\omega^a(t)=\omega\bigl(\inf\{u\ge
t; \pi(u)\in\mathbb{T}^a\}\bigr),
\]
and satisfies $0\le\omega-\omega^a\le a$. Since $\mathbb{T}^a$ is
finite, it
follows from \textup{(\ref{piecewise})} that $\omega^a$ is piecewise
monotone.
Actually, if $U_i^a$ is a time of $[T_i^a,S_i^a]$ at which $\omega$ is
minimal (for $1\le i<N^a$) and if $U_0^a:=0$, $U_{N^a}^a:=1$, then
%
%e2.10 ###
\begin{equation}\label{uia}
\omega(U_i^a)=\omega(S_i^a)-a\qquad \mbox{for $1\le i<N^a$,}
\end{equation}
and
%
%e2.11 ###
\begin{equation}\label{monotone}
\cases{\mbox{$\omega^a$ is nondecreasing on $[U_i^a,T_{i+1}^a]$,}\cr
\mbox{$\omega^a$ is nonincreasing on $[T_i^a,U_i^a]$.}}
\end{equation}

Consider now a general continuous map $\omega$ which does not satisfy
$\pi(0)=\pi(1)$. Then we can again associate the tree $\mathbb{T}$ by
means of
$\delta$ defined by \textup{(\ref{delta})}, but some of the above
properties differ.
However, it is still possible to apply the above discussion to an extended
path $\omega'$ defined on a greater interval, say $[-1,2]$, coinciding with
$\omega$ on $[0,1]$, and satisfying
$\omega'(-1)=\omega'(2)=\inf_{[-1,2]}\omega'$. Then the associated tree
$\mathbb{T}'$ contains $\mathbb{T}$ as a subtree, and the projection
$\pi\dvtx[0,1]\to\mathbb{T}$ is the restriction of $\pi'\dvtx
[-1,2]\to\mathbb{T}'$
to $[0,1]$. Among these paths, we will only consider the minimal extensions;
they are those such that $\mathbb{T}'=\mathbb{T}$. This means that
%
%e2.12 ###
\begin{equation}\label{extension}
\cases{\displaystyle
\omega'(-1)=\omega'(2)=\inf_{[0,1]}\omega,\cr
\mbox{$\omega'$ is nondecreasing on $[-1,0]$, nonincreasing on
$[1,2]$.}}
\end{equation}
Let $U$ be a time of $[0,1]$ at which $\omega$ is minimal and consider
%
%e2.13 ###
\begin{equation}\label{newroot}
O:=\pi(U),\qquad  A:=\pi(0),\qquad  B:=\pi(1)
\end{equation}
(these points are drawn in Figure \ref{fig3} below, in the more general case of paths with jumps).
We choose $O$ as the root of $\mathbb{T}$. Then $O$ belongs to
$[A,B]$, the
points of $[O,A]$ are those such that
$\tau^\nearrow\le0\le\tau^\nwarrow\le1$, and the points of
$[O,B]$ are those
such that $0\le\tau^\nearrow\le1\le\tau^\nwarrow$; for the
points of
$\mathbb{T}\setminus[A,B]$, one has $0<\tau^\nearrow\le\tau
^\nwarrow<1$.

%
%f3 ###
\begin{figure}

\includegraphics{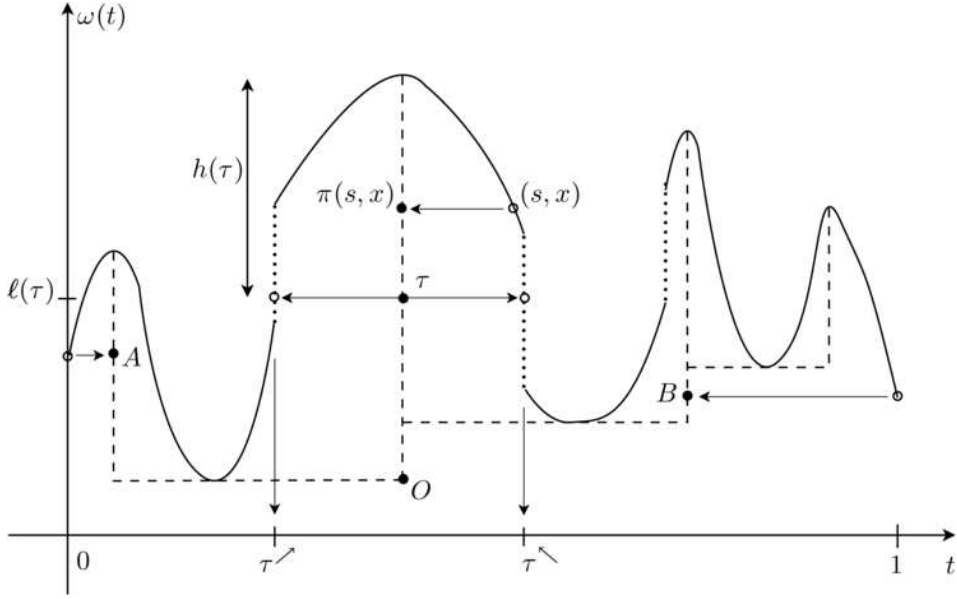}

\caption{A path with jumps, and its tree (dashed lines). The graph
$\mathbb{G}$ is the curve augmented
by the jumps (dotted lines). Are also depicted the map $\pi$ from
$\mathbb{G}$ to $\mathbb{T}$,
the maps $\tau\mapsto\tau^\nearrow$, $\tau\mapsto\tau^\nwarrow$ from
$\mathbb{T}$
to $[0,1]$; in particular, $A=\pi(0,\omega(0))$ and $B=\pi(1,\omega
(1))$.}\label{fig3}
\end{figure}

In particular, if we trim the tree $\mathbb{T}$ and if $\mathbb
{T}^a\ne
\varnothing$,
then the flattened path $\omega^a$ of \textup{(\ref{defomegaa})} is
the restriction
of ${\omega'}^a$ to $[0,1]$. Moreover, the quantities $N^a$, $T_i^a$ and
$S_i^a$ defined in \textup{(\ref{tiasia})} and the similar quantities for
$\omega'$
satisfy
\[
N^a={N'}^a,\qquad  S_i^a=(S')_i^a, \qquad T_i^a=(T')_i^a
\qquad \mbox{for $1\le i<N^a$.}
\]
At $i=N^a$, the time $(T')_{N^a}^a$ may be after time 1, and in this case
$T_{N^a}^a$ is not defined.

%s2.2 ###
\subsection{The length measure on the tree}\label{sec22}

The length measure on $\mathbb{T}$ is the unique measure $\lambda$
which is
supported by the skeleton (the set of leaves have zero measure) and such
that the measure of an arc is equal to its length; in particular, this
measure is $\sigma$-finite and atomless. The existence and uniqueness of
$\lambda$ is elementary for the finite subtrees $\mathbb{T}^a$, and
it is not
difficult to deduce the result for $\mathbb{T}$ by letting
$a\downarrow0$.
It can
be identified to either of the two following measures.

\begin{proposition}\label{lambda12}
Define
\[
\lambda_1:=\int_\mathbb{R}\>\sum_{\tau\in S(\mathbb{T})\dvtx\ell
(\tau
)=x}\delta_\tau
\, dx,\qquad
 \lambda_2:=\int_0^\infty\sum_{\tau\in\partial\mathbb{T}
^a}\delta_\tau \,da
=\int_0^\infty\sum_{\tau\dvtx h(\tau)=a}\delta_\tau\, da,
\]
where $\delta_\tau$ denotes the Dirac mass at $\tau$. Then
$\lambda=\lambda_1=\lambda_2$.
\end{proposition}

Notice that the number of terms in the sum is at most countable for any $x$
in the definition of $\lambda_1$, whereas it is finite for any $a>0$
in the
definition of $\lambda_2$. The integrals are supported by the interval
$[\inf\omega,\sup\omega]$ for the first one, and $[0,\sup\omega
-\inf
\omega]$
for the second one.

\begin{pf*}{Proof of Proposition \protect\ref{lambda12}}
The two measures are supported by the skeleton of the tree; in order to
check that they are equal to $\lambda$, it is sufficient to verify
that they
coincide with it on arcs $[O,\tau]$ for any $\tau$ in the skeleton
$S(\mathbb{T})$. The maps~$\ell$ and $h$ are injective on $[O,\tau
]$, so, if
$\lambda_\mathbb{R}$ denotes the Lebesgue measure on $\mathbb{R}$,
\[
\lambda_1([O,\tau])=\lambda_\mathbb{R}(\ell([O,\tau])
),\qquad
\lambda_2([O,\tau])=\lambda_\mathbb{R}(h([O,\tau])).
\]
Moreover, $\ell$ induces a bijection between $[O,\tau]$ and
$[\ell(O),\ell(\tau)]$, so
\[
\lambda_1([O,\tau])=\ell(\tau)-\ell(O)=\delta(O,\tau)=\lambda
([O,\tau]).
\]
Thus $\lambda_1=\lambda$. For the study of $\lambda_2$, notice that
$h(\tau_0)$ is the distance between $\tau_0$ and any of the highest points
above it. When $\tau_0$ goes from $O$ to $\tau$, then $h(\tau_0)$ is
decreasing; more precisely, it jumps at $\tau_0$, when $\tau_0$ is a
branching point so that no highest point above it is in the direction of
$\tau$; thus $h$ has a finite number of negative jumps, and between these
jumps, it is affine with slope $-1$. Consequently, $h$ induces a bijection
from $[O,\tau]$ onto its image, and this image has Lebesgue measure
$\delta(O,\tau)$. We deduce that $\lambda_2=\lambda$.
\end{pf*}

The measure $\lambda$ is closely related to the two following measures on
$[0,1]$. Say that an excursion begins at time $t$ above level $\omega
(t)$ if
for some $\varepsilon>0$, $\omega(s)>\omega(t)$ for
$t<s<t+\varepsilon$. Let
$E^\nearrow$
be the set of beginnings of excursions above any level; we can define
similarly the set $E^\nwarrow$ of ends of excursions. These two sets
are in
bijection with each other; to each beginning $t$ of an excursion we can
associate its end $\inf\{s>t;\omega(s)=\omega(t)\}$. If we restrict
ourselves to a fixed level $x$, the sets of beginnings and ends of
excursions above $x$ are at most countable, and we can define
%
%e2.14 ###
\begin{equation}\label{omnenw}
\omega^\nearrow:=\int\sum_{s\in
E^\nearrow;\omega(s)=x}\delta_s \,dx,\qquad
\omega^\nwarrow:=\int\sum_{s\in
E^\nwarrow;\omega(s)=x}\delta_s \,dx.
\end{equation}

\begin{proposition}\label{lebesgue}
Assume \textup{(\ref{om01})}. The measures $\omega^\nearrow$ and
$\omega
^\nwarrow$
are $\sigma$-finite and are respectively the images of $\lambda$ by
the maps
$\tau\mapsto\tau^\nearrow$ and $\tau\mapsto\tau^\nwarrow$, and
$\lambda
$ is
the image of $\omega^\nearrow$ and $\omega^\nwarrow$ by the projection
$\pi$. If \textup{(\ref{om01})} does not hold, then, with the notation
\textup{(\ref{newroot})}, the maps $\tau\mapsto\tau^\nearrow$ and
$\tau\mapsto\tau^\nwarrow$ are respectively defined on $\mathbb{T}
\setminus[O,A]$
and $\mathbb{T}\setminus[O,B]$; the relation between $\omega
^\nearrow$ and
$\lambda$ (or between $\omega^\nwarrow$ and $\lambda$) again holds by
restricting $\lambda$ to $\mathbb{T}\setminus[O,A]$ (or $\mathbb
{T}\setminus[O,B]$).
\end{proposition}

\begin{pf}
We only work out the proof under \textup{(\ref{om01})}; the general
case is easily
deduced by considering an extension of $\omega$ satisfying
\textup{(\ref{extension})}. We want to compare the measure $\omega
^\nearrow$ carried
by the set $E^\nearrow$ of beginnings of excursions, with the
measure~$\lambda$ carried by the skeleton $S(\mathbb{T})$. If $s$ is in
$E^\nearrow$,
then $\pi(s)$ is in $S(\mathbb{T})$ and $s=\pi(s)^\nearrow$ except if
$s$ is
at a
local minimum, or the end of a constancy interval of $\omega$; on the other
hand, if $\tau$ is in $S(\mathbb{T})$, then $\tau=\pi(\tau
^\nearrow)$ and
$\tau^\nearrow$ is in $E^\nearrow$ except if it is the beginning of a
constancy interval of $\omega$. Since there are at most countably many local
minima and constancy intervals, we deduce that there exists
$E_0^\nearrow\subset E^\nearrow$ and $S_0(\mathbb{T})\subset
S(\mathbb{T})$
such that
$E^\nearrow\setminus E_0^\nearrow$ and $S(\mathbb{T})\setminus
S_0(\mathbb{T})$
are at
most countable, and the maps $\tau\mapsto\tau^\nearrow$ and $\pi$ are
inverse bijections between $E_0^\nearrow$ and $S_0(\mathbb{T})$. Moreover,
$\lambda$ and $\omega^\nearrow$ are atomless, so they are supported
respectively by $S_0(\mathbb{T})$ and $E_0^\nearrow$. Thus the
relation between
$\lambda$ and $\omega^\nearrow$ claimed in the proposition follows
from this
one-to-one property, the definition \textup{(\ref{omnenw})} of
$\omega
^\nearrow$ and
the property $\lambda=\lambda_1$ of Proposition~\ref{lambda12}. The
case of
$\omega^\nwarrow$ is similar, and the $\sigma$-finiteness follows
from the
$\sigma$-finiteness of $\lambda$.
\end{pf}

We now give a condition on $\mathbb{T}$ with which one can decide whether
$\omega$ has finite or infinite variation (this characterization is also
given in \cite{duquesne07}).

\begin{proposition}\label{finitevar}
The measures $\lambda$, $\omega^\nearrow$ and $\omega^\nwarrow$
are finite
if and only if $\omega$ has finite variation. In this case,
$\omega^\nearrow$ and $\omega^\nwarrow$ are respectively the
positive and
negative parts of the Lebesgue--Stieltjes measure of $\omega$. Moreover,
%
%e2.15 ###
\begin{equation}\label{correction}
\int_0^1|d\omega|=2 \lambda(\mathbb{T})-\delta(0,1).
\end{equation}
\end{proposition}

\begin{pf}
We first work out the proof under the condition \textup{(\ref
{om01})}, so that
$\delta(0,1)=0$. Suppose also that $\mathbb{T}$ is finite, so that
$\lambda$ is
finite and $\omega$ is piecewise monotone [as explained in
\textup{(\ref{piecewise})}]. If, for instance, $\omega$ is nondecreasing on
$[t_1,t_2]$, then it is easily checked from the definitions \textup
{(\ref{omnenw})}
that
\[
\omega^\nearrow([t_1,t_2])=\omega(t_2)-\omega(t_1),\qquad
\omega^\nwarrow([t_1,t_2])=0.
\]
A similar result holds for intervals on which $\omega$ is nonincreasing, so
we deduce that the proposition holds true in this case. If $\mathbb
{T}$ is not
finite, consider the tree $\mathbb{T}^a$ of~\textup{(\ref
{deftreea})} and its path
$\omega^a$ of \textup{(\ref{defomegaa})}. Notice that
$\lambda(\mathbb{T}^a)\uparrow\lambda(\mathbb{T})$ as $a\downarrow
0$. For
$b<a$, one
has $\mathbb{T}^a\subset\mathbb{T}^b$, and the path $\omega^a$ is
obtained from
$\omega^b$ by a change of time, so the variation of $\omega^a$
increases as
$a$ decreases, and is bounded by the variation of $\omega$; since the
variation is a lower semicontinuous function of the path, it follows that
the variation of $\omega^a$ converges to the variation of $\omega$ as
$a\downarrow0$, so
%
%e2.16 ###
\begin{equation}\label{total}
\int_0^1|d\omega|=\lim\int_0^1|d\omega^a|=2\lim\lambda(\mathbb{T}
^a)=2\lambda(\mathbb{T})
\end{equation}
(we have applied the first part of the proof to $\omega^a$ and
$\mathbb{T}^a$).
Thus $\omega$ has finite variation if and only if $\lambda$ is finite.
Moreover, if $\omega$ has finite variation, one checks similarly that the
positive part $d\omega^+$ of the Lebesgue--Stieltjes measure of $\omega$
satisfies
\[
(d\omega^+)([s,t])=\lim(d\omega^a)^+([s,t])
=\lim\omega^\nearrow\bigl([s,t]\cap\pi^{-1}(\mathbb{T}^a)\bigr)
=\omega^\nearrow([s,t])
\]
where we have used the fact that $(\omega^a)^\nearrow$ is the
restriction of
$\omega^\nearrow$ to $\pi^{-1}(\mathbb{T}^a)$. If~\textup{(\ref
{om01})} does not
hold, we
can consider an extension of $\omega$ satisfying \textup{(\ref{extension})}
and then
restrict to $[0,1]$. In this case, with the notation \textup{(\ref
{newroot})}, the
points $\tau$ of $[A,B]$ are such that $\tau^\nearrow\le0$ or
$\tau^\nwarrow\le1$ and should not be counted twice in the total variation
of~$\omega$ in~\textup{(\ref{total})}. The correction which has to
be made is
$\lambda([A,B])=\delta(0,1)$, so we obtain~\textup{(\ref{correction})}.
\end{pf}

%s2.3 ###
\subsection{Paths with jumps}\label{sec23}

Let us explain how our construction of $\mathbb{T}$ can be extended to c\`{a}dl\`{a}g
paths $\omega$ (paths which are right-continuous and have left
limits), see Figure \ref{fig3}; we
apply the classical idea of embedding these paths into continuous paths by
opening temporal windows at times of jumps and considering interpolated
continuous paths (this idea has been used for the rough paths theory in
\cite{williams01}).

Let $\mathbb{G}$ be the set of points $(t,x)$ such that $0\le t\le1$ and
$x$ is
between $\omega(t-)$ and~$\omega(t)$. This is the graph of $\omega$
augmented by the segments joining $(t,\omega(t-))$ and $(t,\omega
(t))$. Then
define
\[
\delta((t,x),(t,x')):=|x'-x|
\]
and
\[
\delta((s,x),(t,x')):=x+x'-2\biggl(\inf_{(s,t)}\omega
\wedge
x\wedge x'\biggr)
\]
if $s<t$. If $\omega$ is continuous, then $\mathbb{G}$ and $[0,1]$
are naturally
identified, in such a way that $\delta$ coincides with the previous
definition \textup{(\ref{delta})}.

Let us say that two points of $\mathbb{G}$ satisfy $(t,x)\le(t',x')$
if either
$t<t'$, or $t=t'$ and $x$ is between $\omega(t-)$ and $x'$. This is a total
order, and $\mathbb{G}$ can be endowed with the topology generated by open
intervals for this order; actually, this topology coincides with the
topology of $\mathbb{G}$ considered as a subset of $\mathbb{R}^2$.

\begin{proposition}\label{graphtree}
The map $\delta$ is a semi-distance on $\mathbb{G}$, and
$\mathbb{T}=(\mathbb{G}/\delta,\delta)$ is a compact real tree.
Actually, there
exists a continuous map $\omega'$ such that $\omega$ is obtained from
$\omega'$ by an increasing (not necessarily surjective) time change, and $\mathbb{T}$
is the
tree associated to $\omega'$.
\end{proposition}

\begin{pf}
Suppose that $\omega$ is not continuous (the result is evident otherwise).
Let~$J$ be the set of times where $\omega$ jumps, and let $(S(t);t\in
J)$ be
a family of (strictly) positive numbers such that $\sum S(t)=1$. Let
\[
\Lambda(t,x):=\dfrac12\Biggl(t+\sum_{u<t}S(u)+S(t)
\frac{x-\omega(t-)}{\omega(t)-\omega(t-)}1_J(t)\Biggr).
\]
Then $\Lambda$ is an increasing bijection from $\mathbb{G}$ onto
$[0,1]$, so
$\mathbb{G}$ and $[0,1]$ can be identified, and previous results on
the tree
representation for continuous functions defined on $[0,1]$ can also be
applied to continuous functions on $\mathbb{G}$. Thus, in order to
prove the
proposition, it is sufficient to find a map $\omega'$ defined on
$\mathbb{G}$.
Put $\omega'(t,x):=x$. It induces the semi-distance
\[
\omega'(s,x)+\omega'(t,x')-2\inf_{[(s,x),(t,x')]}\omega'
=\delta((s,x),(t,x')),
\]
so its tree is $\mathbb{T}$. Moreover, $\omega=\omega'\circ Q$ for
the increasing
time change $Q(t):=(t,\omega(t))$.
\end{pf}

In this setting, let $\pi$ be the projection of $\mathbb{G}$ on
$\mathbb{T}$. We
extend the notation \textup{(\ref{newroot})} by
\[
O:=\pi\bigl(U,\omega(U)\wedge\omega(U-)\bigr),\qquad  A:=\pi(0,\omega
(0)),\qquad
B:=\pi(1,\omega(1)),
\]
where $U$ is a time at which $\omega(U)\wedge\omega(U-)=\inf\omega
$. Let
$E^\nearrow$ be the set of $(t,x)$ in~$\mathbb{G}$ such that $\omega
(s)>x$ for
any $t<s<t+\varepsilon$ and some $\varepsilon>0$, define $E^\nwarrow
$ similarly,
and let
\[
\omega^\nearrow:=\int\sum_{s\dvtx(s,x)\in E^\nearrow}\delta_s\,
dx,\qquad
\omega^\nwarrow:=\int\sum_{s\dvtx(s,x)\in E^\nwarrow}\delta_s\, dx.
\]
Notice also that all the points of $\pi^{-1}(\tau)$ are at the same level;
we let $\tau^\nearrow$ and $\tau^\nwarrow$ be the infimum and
supremum of
the time component of this set.

\begin{proposition}
The measures $\omega^\nearrow$ and $\omega^\nwarrow$ are the images of
$\lambda$ by $\tau\mapsto\tau^\nearrow$ and $\tau\mapsto\tau
^\nwarrow$
[after restricting $\lambda$ as in Proposition \ref{lebesgue} if
\textup{(\ref{om01})} does not hold]. The statements of Proposition
\ref{finitevar}
about the finite variation case again hold true.
\end{proposition}

\begin{pf}
Let us use the notation of the proof of Proposition \ref{graphtree}.
The set
$E^\nearrow$ is the set of beginnings of excursions of $\omega'$, so
$\omega^\nearrow$ is the projection on the time component of
$(\omega')^\nearrow$; we deduce the first statement. Moreover,
$\omega
'$ is
monotone on the intervals corresponding to the jumps of $Q$, so the total
variations of $\omega$ and $\omega'$ coincide (a more general result
will be
proved in Theorem \ref{vpomp}), and the Lebesgue--Stieltjes measure of
$\omega$ is again deduced from its analogue for $\omega'$ by
projection on
the time component.
\end{pf}

%s3 ###
\section{$p$-variation and trees}
\label{pvariation}

Let us now assume that $\omega$ has finite $p$-variation for some
$p\ge1$,
so that
%
%e3.1 ###
\begin{equation}\label{pvarfinite}
V_p(\omega):=\sup_{(t_i)}V_p(\omega,(t_i))
:=\sup_{(t_i)}\sum_i|\omega(t_{i+1})-\omega(t_i)
|^p<\infty,
\end{equation}
where the supremum is with respect to all the subdivisions of $[0,1]$
(notice that a nonconstant continuous map cannot have finite $p$-variation
for $p<1$). Let us first assume that $\omega$ is continuous (the
c\`{a}dl\`{a}g case
will be dealt with in Section \ref{jumps}). We first want to
describe the
property \textup{(\ref{pvarfinite})} by means of the geometry of
$\mathbb{T}$. In
particular, $V_p(\omega)<\infty$ implies $V_q(\omega)<\infty$ for
$q\ge p$,
and we are interested in the variation index
%
%e3.2 ###
\begin{equation}\label{varindex}
\mathcal{V}(\omega):=\inf\{p\ge1;V_p(\omega)<\infty\}.
\end{equation}

%s3.1 ###
\subsection{The variation index}\label{sec31}

Let us recall that we have defined in \textup{(\ref{deftreea})} approximations
$\mathbb{T}^a$ of $\mathbb{T}$ obtained by trimming the tree, that $N^a$,
defined by
\textup{(\ref{tiasia})}, is the number of leaves of $\mathbb{T}^a$,
and that the flattened
path $\omega^a$ of \textup{(\ref{defomegaa})} is associated to
$\mathbb{T}^a$; let
$L^a:=\lambda(\mathbb{T}^a)$ be its total length. As $a\downarrow0$,
each branch
of $\mathbb{T}^a$ grows at unit speed at its leaves, so
%
%e3.3 ###
\begin{equation}\label{lanb}
L^a=\int_a^\infty N^b\,db.
\end{equation}
If $\pi(0)=\pi(1)$, we deduce from Proposition \ref{finitevar} that
$L^a$ is
the mass of the positive part of $d\omega^a$, so, by applying
\textup{(\ref{monotone})} and \textup{(\ref{uia})},
\begin{eqnarray*}
L^a&=&\omega(T_1^a)+\sum_{i=2}^{N^a}\bigl(\omega(T_i^a)-\omega
(U_{i-1}^a)\bigr)\\
&=&\omega(T_1^a)+\sum_{i=2}^{N^a}\bigl(\omega(T_i^a)-\omega
(S_{i-1}^a)+a\bigr)\\
&=&\sum_{i=1}^{N^a}\bigl(\omega(T_i^a)-\omega(S_{i-1}^a)\bigr)+(N^a-1)a.
\end{eqnarray*}
If $\pi(0)\ne\pi(1)$, then this equation has to be corrected as in
\textup{(\ref{correction})}; notice, however, that the correction is
bounded, so if
$\omega$ has infinite variation, then
%
%e3.4 ###
\begin{equation}\label{lasum}
L^a\sim\sum_{i=1}^{N^a-1}\bigl(\omega(T_i^a)-\omega
(S_{i-1}^a)
\bigr)+a N^a
\qquad \mbox{as $a\downarrow0$.}
\end{equation}
Thus $L^a$ is easily estimated from the path $\omega$, the times
$S_i^a$ and
$T_i^a$, and the number~$N^a$ of \textup{(\ref{tiasia})}.

We consider two other metric characteristics of $\mathbb{T}$, namely
its upper
box (or Minkowski) dimension (see, for instance, \cite{falconer90})
defined by
\[
\operatorname{\overline{dim}}\mathbb{T}:=\limsup_{a\downarrow
0}\frac{\log\mathcal{N}(a)}{\log(1/a)}
\]
where $\mathcal{N}(a)$ is the minimal number of balls of radius $a$
which are
needed to cover~$\mathbb{T}$, and the index
\[
\mathcal{H}(\mathbb{T}):=\inf\biggl\{p\ge1;
\int_\mathbb{T}(h(\tau))^{p-1}\lambda(d\tau)<\infty
\biggr\}
\]
where $h(\tau)$ is the height of the highest branch above $\tau$. The
aim of
this subsection is to prove that all these quantities are related to the
variation index $\mathcal{V}(\omega)$ defined in \textup{(\ref
{varindex})}, and in
particular prove the result announced in \textup{(\ref{vpomega})}.

\begin{theorem}\label{dimtree}
Let $\omega$ be a (nonconstant) continuous function. Then
\begin{eqnarray*}
\mathcal{V}(\omega)
&=&\mathcal{H}(\mathbb{T})\\
&=&\limsup_{a\to0}\frac{\log L^a}{\log(1/a)}+1\\
&=&\limsup_{a\to0}\frac{\log N^a}{\log(1/a)}\vee1\\
&=&\operatorname{\overline{dim}}\mathbb{T}.
\end{eqnarray*}
\end{theorem}

\begin{pf}
Denoting by $I_1,\ldots,I_5$ the successive terms of the theorem, we prove
that
\[
I_1\le I_2\le I_3\le I_4\le I_5\le I_1.
\]
These five inequalities are proved in the five following steps.

\textsc{Proof of $I_1\le I_2$.} Let $s<t$ be two times, and let
$\tau_0$ be the most recent common ancestor of $\pi(s)$ and $\pi
(t)$. Then
\[
\ell(\tau_0)=\min_{[s,t]}\omega
\]
so
\[
|\omega(t)-\omega(s)|
\le\max\bigl(\omega(s)-\ell(\tau_0),\omega(t)-\ell(\tau
_0)\bigr)
\]
and
\[
|\omega(t)-\omega(s)|^p\le\bigl(\omega(s)-\ell(\tau
_0)\bigr)^p
+\bigl(\omega(t)-\ell(\tau_0)\bigr)^p.
\]
On the other hand,
\begin{eqnarray*}
\bigl(\omega(t)-\ell(\tau_0)\bigr)^p
&=&p\int_{[\tau_0,\pi(t)]}\bigl(\omega(t)-\ell(\tau)
\bigr)^{p-1}\lambda
(d\tau)\\
&\le& p\int_{[\tau_0,\pi(t)]}\bigl(h(\tau)-h(\pi(t))
\bigr)^{p-1}\lambda
(d\tau)\\
&\le& p\int_{[\tau_0,\pi(t)]}(h(\tau))^{p-1}\lambda
(d\tau)
\end{eqnarray*}
where we have used in the second line the property
\begin{eqnarray*}
h(\tau)-h(\pi(t))&=&\max_{[\tau^\nearrow,\tau^\nwarrow]}\ell
-\ell(\tau)
-\max_{[\pi(t)^\nearrow,\pi(t)^\nwarrow]}\ell+\omega(t)\\
&\ge&\omega(t)-\ell(\tau)
\end{eqnarray*}
valid for $\pi(t)$ above $\tau$. The same property holds at time $s$,
so by
addition,
\[
|\omega(t)-\omega(s)|^p\le p\int_{[\pi(s),\pi(t)]}
(h(\tau))^{p-1}\lambda(d\tau).
\]
If $(t_i)$ is a subdivision of $[0,1]$, we can sum up these estimates for
$s=t_i$ and $t=t_{i+1}$. Since almost any $\tau$ appears at most twice in
the right-hand sides (at times $\tau^\nearrow$ and $\tau^\nwarrow
$), we
deduce
%
%e3.5 ###
\begin{equation}\label{vphtau}
V_p(\omega)\le2p\int_\mathbb{T}(h(\tau))^{p-1}\lambda
(d\tau).
\end{equation}
In particular $I_1\le I_2$.

\textsc{Proof of $I_2\le I_3$.} It follows from $\lambda
=\lambda_2$
(Proposition \ref{lambda12}) and from \textup{(\ref{lanb})} that for $p>1$,
\[
\int_\mathbb{T}(h(\tau))^{p-1}\lambda(d\tau)=\int
_0^\infty
a^{p-1}N^a\,da
=(p-1)\int_0^\infty a^{p-2}L^a\,da.
\]
We deduce that if $L^a\le Ca^{1-\kappa}$ for some $\kappa<p$, then the
integral is finite, so \mbox{$I_2\le I_3$}.

\textsc{Proof of $I_3\le I_4$.} This inequality follows from
\textup{(\ref{lanb})}.

\textsc{Proof of $I_4\le I_5$.} Above each $\tau\in
\partial
\mathbb{T}^a$
there is a $\tau'$ such that $\delta(\tau,\tau')=a$, and the $N^a$ balls
with centers $\tau'$ and radius $a$ are disjoint; this implies that the
number of balls of radius $a/2$ which is needed to cover $\mathbb{T}$
is at least
$N^a$; we also have $\operatorname{\overline{dim}}\mathbb{T}\ge1$,
so we deduce that $I_4\le I_5$.

\textsc{Proof of $I_5\le I_1$.} For $a>0$, let $t_0=0$ and
\[
t_{i+1}=\inf\{t\ge t_i;|\omega(t)-\omega(t_i)|\ge a\}.
\]
Let $\tau_i$ be the most recent common ancestor of $\pi(t_i)$ and
$\pi(t_{i+1})$, so that $\ell(\tau_i)=\inf_{[t_i,t_{i+1}]}\omega$. Consider
the closed ball $B_i$ of $\mathbb{T}$ with center $\tau_i$ and with
radius $2a$,
so that $\pi([t_i,t_{i+1}])$ is included in this ball. Then the union of
$B_i$ is a covering of~$\mathbb{T}$. Moreover, the number of these
balls is
dominated by $V_p(\omega)/a^p$, so the upper box dimension of $\mathbb
{T}$ is
dominated by $p$ as soon as $p>\mathcal{V}(\omega)$. We deduce that
$I_5\le I_1$.
\end{pf}

\begin{remark}
If $\pi(0)=\pi(1)$, we have
\[
V_p(\omega)\ge\sum_{\tau\in\partial\mathbb{T}^a}
\biggl(\biggl(\sup_{[\tau^\nearrow,\tau^\nwarrow]}\omega-\omega
(\tau
^\nearrow)\biggr)^p
+\biggl(\sup_{[\tau^\nearrow,\tau^\nwarrow]}\omega-\omega(\tau
^\nwarrow
)\biggr)^p\biggr)
=2a^pN^a.
\]
If $\pi(0)\ne\pi(1)$, we have to omit the first term for the first
leaf of
$\mathbb{T}^a$ ($\tau^\nearrow$ may be before time 0), and the second
term for
the last leaf of $\mathbb{T}^a$ ($\tau^\nwarrow$ may be after time~1). Thus
%
%e3.6 ###
\begin{equation}\label{vpapna}
V_p(\omega)\ge a^pN^a
\end{equation}
and the right-hand side can be doubled if $\pi(0)=\pi(1)$.
\end{remark}

\begin{remark}
Other related estimates of $V_p(\omega)$ using numbers of upcrossings were
previously known; see \cite{bruneau79,stricker79}.
\end{remark}

\begin{remark}
The link between the dimension of $\mathbb{T}$ and the behavior of
$N^a$ is
similar to the link between the dimension of the boundary of discrete trees
and their growth (see page 201 of \cite{peres99}).
\end{remark}

\begin{remark}
A more classical fractal dimension related to a path $\omega$ is the
dimension of its graph as a subset of $\mathbb{R}^2$. This dimension
(which is
bounded by 2) is of course generally different from the dimension of
$\mathbb{T}$.
\end{remark}

Other well-known notions of dimensions (\cite{falconer90}) are the packing
dimension $\dim_P\mathbb{T}$ and the Hausdorff dimension $\dim
_H\mathbb{T}$,
and we
always have
%
%e3.7 ###
\begin{equation}\label{dimhdimp}
\dim_H\mathbb{T}\le\dim_P\mathbb{T}\le\operatorname{\overline{dim}}\mathbb{T}.
\end{equation}
Some of these inequalities may be strict. For instance, consider the path
$\omega$ which is affine on each interval $[1/(n+1),1/n]$, and such that
\[
\omega\bigl(1/(2k+1)\bigr)=0,\qquad \omega\bigl(1/(2k)\bigr)=1/k^\alpha.
\]
Then
\[
V_p(\omega)=2\sum k^{-\alpha p}
\]
for $p\ge1$, so $\operatorname{\overline{dim}}\mathbb{T}=\mathcal
{V}(\omega)=1/\alpha\vee1$.
On the other
hand, the tree is a star with a countable number of branches, and its
Hausdorff and packing dimensions are 1.

On the other hand, if $\omega$ has the same variation index
$\mathcal{V}(\omega)$ on any interval $[s,t]$ with $s<t$, then any
open subset
of $\mathbb{T}$ has the same upper box dimension, so in this case
(\cite{falconer90})
\[
\dim_P\mathbb{T}=\operatorname{\overline{dim}}\mathbb{T}=\mathcal
{V}(\omega).
\]

We will now see an example where the Hausdorff dimension is also equal
to~$\mathcal{V}(\omega)$.

%s3.2 ###
\subsection{The fractional Brownian case}\label{sec32}

We now consider the case where $\omega$ is a typical path of a fractional
Brownian motion $W$. This is a centered Gaussian process $(W_t;t\in
\mathbb{R})$
with covariance function
\[
\operatorname{cov}(W_s,W_t)=\frac{\sigma^2}{2}
(|s|^{2H}+|t|^{2H}-|t-s|^{2H})
\]
for the Hurst parameter $0<H<1$ and the coefficient $\sigma^2>0$. It
satisfies the scaling property
%
%e3.8 ###
\begin{equation}\label{scaling}
(W_{ct}; t\in\mathbb{R})\stackrel{\mathrm{law}}{=}(c^HW_t; t\in
\mathbb{R})
\end{equation}
for $c>0$. In this subsection, we let $\omega$ be a path of $W$ restricted
to $[0,1]$ and extended to $[-1,2]$ by the technique of \textup{(\ref
{extension})};
we compute the Hausdorff dimension of the tree $\mathbb{T}$, and
describe an
estimator of $H$ based on $\mathbb{T}$.

The property $V_p(W)<\infty$ for $p>1/H$ is well known; it is classically
obtained from the $(1/p)$-H\"{o}lder continuity of the paths, which itself is
obtained by means of the Kolmogorov criterion and the estimation
\[
\|W_t-W_s\|_q=C_q\sigma(t-s)^H
\]
on the $L^q$ norm of the increments for any $q\ge1$. It actually follows
from this estimation that the moments of $V_p(W)$ are finite.

\begin{proposition}\label{dimfbm}
For almost any path $\omega$ of $W$, one has
\[
\dim_H\mathbb{T}=\mathcal{V}(\omega)=1/H.
\]
\end{proposition}

\begin{pf}
The property $\mathcal{V}(W)\le1/H$ follows from the discussion
preceding the
proposition. From \textup{(\ref{dimhdimp})} and Theorem \ref
{dimtree}, it is
therefore sufficient to prove that $\dim_H\mathbb{T}\ge1/H$. The constants
involved in this proof depend on $H$ and $\sigma$. It is known from
\cite{molchan99} that
%
%e3.9 ###
\begin{equation}\label{pinf}
\mathbb{P}\biggl[\inf_{[0,1/2]}W>-u\biggr]=O(u^\gamma)
\end{equation}
as $u\downarrow0$, for any $\gamma<1/H-1$. Moreover, if $(\mathcal{F}
_t;0\le
t\le1)$ is the filtration of $W$, the conditional law of $W_1-W_{1/2}$ given
$\mathcal{F}_{1/2}$ is a Gaussian law with deterministic positive
variance, so
%
%e3.10 ###
\begin{equation}\label{pw1}
\mathbb{P}[|W_1|<u\mid\mathcal{F}_{1/2}]\le C u.
\end{equation}
The event $\{\delta(0,1)<u\}$ is included in the intersection of the two
events of \textup{(\ref{pinf})} and \textup{(\ref{pw1})}, so
\[
\mathbb{P}[\delta(0,1)<u]=O(u^{\gamma+1}).
\]
We deduce that $\delta(0,1)^{-p}$ is integrable for $p<1/H$. From the
scaling property~\textup{(\ref{scaling})}, $\delta(s,t)^{-p}$ is also
integrable, and
\[
\mathbb{E}\delta(s,t)^{-p}=C(t-s)^{-pH},
\]
so
%
%e3.11 ###
\begin{equation}\label{energy}\qquad
\mathbb{E}\int\!\!\int_{\mathbb{T}\times\mathbb{T}}\delta(\tau
_1,\tau
_2)^{-p}\nu(d\tau_1) \nu
(d\tau_2)
=\mathbb{E}\int_0^1\int_0^1\delta(s,t)^{-p}\,ds\,dt<\infty
\end{equation}
for the projection $\nu$ of the Lebesgue measure of $[0,1]$ on
$\mathbb{T}
$. The
double integral of the left-hand side is the $p$-energy of the measure
$\nu$
on the metric space $(\mathbb{T},\delta)$. Its almost sure finiteness implies
that $\dim_H\mathbb{T}\ge p$ for any $p<1/H$ (see, for
instance,~\cite{falconer90}), so $\dim_H\mathbb{T}\ge1/H$.
\end{pf}

Dimensions of L\'{e}vy trees have been computed in \cite{duqleg05}. This
includes our tree~$\mathbb{T}$ for $H=1/2$, and for this tree, the exact
Hausdorff measure has been obtained in~\cite{duqleg07}. Here, we do not
look for a so precise result, but verify that the normalization of the
length measure $\lambda$ on $\mathbb{T}^a$ converges to the measure
$\nu$ of the
previous proof; the same property is verified for the uniform measure on
leaves of $\mathbb{T}^a$. In this sense, $\nu$ can be viewed as a
uniform measure
on the leaves of the tree. This will be a corollary of the following result
(Proposition \ref{nu12}).

\begin{proposition}\label{nala}
For almost any path $\omega$ of $W$, we have
\[
N^a\sim C(H) \sigma^{1/H}a^{-1/H},\qquad  L^a\sim C(H) \frac
{H}{1-H}\sigma^{1/H}a^{1-1/H}
\]
as $a\downarrow0$, for some $C(H)>0$.
\end{proposition}

\begin{pf}
Since $N^a$ and $L^a$ are related to each other by means of \textup
{(\ref{lanb})}, it
is sufficient to study $N^a$. Moreover, $\sigma$ acts as a multiplicative
coefficient on the path, so $N^a$ for the process with parameter
$\sigma$
has the same law as $N^{a/\sigma}$ for the process with parameter 1; thus
it is sufficient to consider the case $\sigma=1$. If $p>1/H$, it follows
from the finiteness of the moments of $V_p(W)$ and from \textup{(\ref
{vpapna})} that
%
%e3.12 ###
\begin{equation}\label{naq}
\|N^a\|_q\le C a^{-p}
\end{equation}
for any $q\ge1$ and some $C=C(p,q,H)$. In the two following steps, we study
successively the expectation and the variance of $N^a$.

\textsc{Study of $\mathbb{E}[N^a]$.} Consider in this proof
the whole
path $(\omega(t);t\in\mathbb{R})$ of $W$, and its associated (noncompact) tree
$\mathbb{T}_{-\infty,+\infty}$. For $s<t$, let $N_{s,t}^a$, respectively
$\widetilde N_{s,t}^a$, be the numbers of leaves of the trimmed tree
$\mathbb{T}_{-\infty,+\infty}^a$ such that $s<\tau^\nearrow<\tau
^\nwarrow<t$,
respectively $s\le\tau^\nearrow<t$. Then
%
%e3.13 ###
\begin{equation}\label{ntildea}
\widetilde N_{s,t}^a-N_{s,t}^a\in\{0,1\},\qquad
N^a-N_{0,1}^a\in\{0,1,2\}
\end{equation}
(actually one may have $\widetilde N_{s,t}^a-N_{s,t}^a=2$ if $s$ is some
$\tau^\nearrow$, but this happens with zero probability for any fixed $s$).
On the other hand, it follows from the scaling property~\textup{(\ref
{scaling})} of
$W$ that
\[
\mathbb{E}[\widetilde N_{0,1}^a]=\mathbb{E}[\widetilde N_{0,a^{-1/H}}^1].
\]
The law of $W$ is shift invariant and $[s,t)\mapsto\widetilde
N_{s,t}^1$ is
additive, so $\mathbb{E}[\widetilde N_{s,t}^1]$ is proportional to
$t-s$, and
\[
\mathbb{E}[\widetilde N_{0,1}^a]=a^{-1/H}\mathbb{E}[\widetilde N_{0,1}^1].
\]
Thus the result of the proposition holds in expectation for
$C(H)=\mathbb{E}[\widetilde N_{0,1}^1]$.

\textsc{Study of $\operatorname{var}(N^a)$.} It follows
from \textup{(\ref{ntildea})} and
the additivity of $[s,t)\mapsto\widetilde N_{s,t}^a$ that
\[
|N_{s,u}^a+N_{u,t}^a-N_{s,t}^a|\le2
\]
for $s\le u\le t$. Thus, by considering a regular subdivision of $[0,1]$
with mesh $\Delta t$, we have
%
%e3.14 ###
\begin{equation}\label{nana}
\Bigl|N^a-\sum N_{t_i,t_{i+1}}^a\Bigr|\le2 \Delta t^{-1}+2.
\end{equation}
Moreover,
\begin{eqnarray*}
\operatorname{var}(N_{t_i,t_{i+1}}^a)
&=&\operatorname
{var}(N_{0,\Delta
t}^a)
=\operatorname{var}(N_{0,1}^{a \Delta t^{-H}})\\
&\le&\mathbb{E}[(N_{0,1}^{a \Delta t^{-H}})^2]
\le\mathbb{E}[(N^{a \Delta t^{-H}})^2]
\le C a^{-2p}(\Delta t)^{2pH}
\end{eqnarray*}
for $p>1/H$, where we have used the scaling property in the second equality,
and~\textup{(\ref{naq})} in the last inequality. Since
$N_{t_i,t_{i+1}}^a$ depends
only on the increments of~$\omega$ on $[t_i,t_{i+1}]$, we deduce from the
result \textup{(\ref{rhofjfk})} of Appendix \ref{mixing} that
%
%e3.15 ###
\begin{eqnarray}\label{varnti}
\operatorname{var}\Bigl(\sum N_{t_i,t_{i+1}}^a\Bigr)
&\le& C a^{-2p}(\Delta t)^{2pH}\sum_{k,j\le\Delta t^{-1}}\frac
1{1+|k-j|^{1-H}}\nonumber\\[-8pt]\\[-8pt]
&\le& C' a^{-2p}(\Delta t)^{2pH-H-1},\nonumber
\end{eqnarray}
so, by joining \textup{(\ref{nana})} and \textup{(\ref{varnti})},
\[
\operatorname{var}(N^a)\le C\bigl(a^{-2p}(\Delta t)^{2pH-H-1}+(\Delta
t)^{-2}\bigr).
\]
We choose $\Delta t\sim a^\alpha$ for $0<\alpha<1/H$, so
\[
\operatorname{var}(N^a)\le C\bigl(a^{-2\alpha}+a^{-2p+\alpha
(2pH-H-1)}\bigr).
\]
By choosing $p$ and $\alpha$ close enough to $1/H$, we have
\[
\operatorname{var}(N^a)\le C a^{2\varepsilon-2/H}
\]
for some $\varepsilon>0$.

\textsc{Conclusion of the proof.} The two previous steps
show that
$a^{1/H}N^a$ converges in $L^2$ to a constant, and that the rate of
convergence is at most of order~$a^\varepsilon$. From the Borel--Cantelli
lemma, the
convergence is almost sure on a sequence $a_n=n^{-\beta}$ for $\beta$ large
enough. Since $a\mapsto N^a$ is monotone, we deduce from
\[
a_{n+1}^{1/H}N^{a_n}\le a^{1/H}N^a\le a_n^{1/H}N^{a_{n+1}}
\]
for $a_{n+1}\le a\le a_n$, that the convergence is actually almost sure as
$a\downarrow0$.
\end{pf}

\begin{proposition}\label{nu12}
For almost any path $\omega$ of $W$, the measures
\[
\nu_1^a:=\frac1{N^a}\sum_{\tau\in\partial\mathbb{T}^a}\delta
_\tau
 \quad\mbox{and}\quad
\nu_2^a:=\frac1{L^a} \lambda|_{\mathbb{T}^a}
\]
converge weakly to the projection $\nu$ on $\mathbb{T}$ of the
Lebesgue measure
of $[0,1]$.
\end{proposition}

\begin{pf}
Let $\mu_1^a$ and $\mu_2^a$ be the images of $\nu_1^a$ and $\nu
_2^a$ by
$\tau\mapsto\tau^\nearrow$. One has $\pi(\tau^\nearrow)=\tau$,
so $\nu_1^a$
and $\nu_2^a$ are the images of $\mu_1^a$ and $\mu_2^a$ by $\pi$. Since
$\pi$ is continuous, it is sufficient to prove that $\mu_1^a$ and
$\mu_2^a$
converge weakly to the Lebesgue measure of $[0,1]$, and therefore that
$\mu_1^a([s,t])$ and $\mu_2^a([s,t])$ converge to $t-s$. But
$\mu_1^a([s,t])$ counts the proportion of leaves of $\mathbb{T}^a$
which satisfy
$s\le\tau^\nearrow\le t$; the number of such leaves is close to the number
$\widetilde N_{s,t}^a$ of the proof of Proposition \ref{nala}; it can be
estimated from Proposition \ref{nala} and the scaling property, and we can
conclude. The study of $\mu_2^a$ is similar.
\end{pf}

We can deduce estimators for $H$ from Proposition \ref{nala}. Our
result is
an alternative to the generalized quadratic variation approach
\cite{istaslang97}. For instance, we can consider $N^{2a}/N^a$ or
$L^{2a}/L^a$, so that the unknown coefficient $\sigma$ is eliminated.
However, we can also use
\[
\lim_{a\downarrow0}\frac{a N^a}{L^a}=\frac1H-1.
\]
Roughly speaking, the estimator $a N^a/L^a$ counts the normalized
number of
changes in the sense of variation of $\omega^a$. The smaller $H$ is, the
more often the sense of variation of $\omega^a$ changes. From \textup
{(\ref{lasum})},
we deduce the following result.

\begin{proposition}\label{estimator}
The Hurst parameter $H$ of the fractional Brownian motion $(W_t; 0\le
t\le1)$ can be estimated from the relation
\[
\lim_{a\downarrow0}\frac{1}{a N^a}\sum_{i=1}^{N^a-1}
(W_{T_i^a}-W_{S_{i-1}^a})=\frac{2H-1}{1-H}
\]
which holds almost surely, where $N^a$, $S_i^a$, $T_i^a$ were defined in
\textup{(\ref{tiasia})}.
\end{proposition}

%s3.3 ###
\subsection{The case with jumps}\label{jumps}

We now consider a c\`{a}dl\`{a}g path $\omega$. We have seen in Proposition
\ref{graphtree} how it can be written as a time-changed path
$\omega=\omega'\circ Q$ for a continuous $\omega'$ defined on
$\mathbb{G}$, and
the trees of $\omega$ and $\omega'$ coincide. Actually, the
variations also
coincide, so the tree $\mathbb{T}$ can again be used to study the
variations of~$\omega$.

\begin{theorem}\label{vpomp}
Let $\omega$ be a c\`{a}dl\`{a}g path and $\omega'$ the associated continuous path.
One has $V_p(\omega)=V_p(\omega')$ for any $p\ge1$. In particular,
$\mathcal{V}(\omega')=\mathcal{V}(\omega)$ and Theorem \ref
{dimtree} again
holds.
\end{theorem}

\begin{pf}
The relation $\omega=\omega'\circ Q$ immediately implies $V_p(\omega
)\le
V_p(\omega')$. In order to verify the reverse inequality, we notice that
when computing $V_p(\omega')$, it is sufficient to consider subdivisions
$(t_i)$ consisting of local extrema of $\omega'$; thus these times are in
the closure of the image of $Q$; consequently, from the continuity of
$\omega'$, it is sufficient to consider times in the image of $Q$, so that
we can conclude.
\end{pf}

We now give applications of the tree representation to martingales and L\'{e}vy
processes. In the following result, we recover with our method a result of
\cite{pisierxu88} (which was given in discrete time). Notice, however, that
our results are only for the real-valued case, whereas \cite{pisierxu88}
considers the Banach space-valued case.

\begin{proposition}\label{super}
Consider a purely discontinuous martingale $X=(X_t; 0\le t\le1)$ for a
filtration $(\mathcal{F}_t; 0\le t\le1)$. Let $1<p<2$; then
%
%e3.16 ###
\begin{equation}\label{evp}
\mathbb{E}[V_p(X)]\le C_p \mathbb{E}\sum|\Delta X_t|^p.
\end{equation}
\end{proposition}

\begin{pf}
The proof is divided into two steps; in the first step, we reduce the
problem to a particular case.

\begin{step}\label{step1}
Let $S_0:=0$ and $(S_k; k\ge1)$ be the times of
jumps of an independent standard Poisson process, and consider
\[
X_t^\varepsilon=\sum X_{\varepsilon S_k}1_{\{\varepsilon S_k\le
t<\varepsilon S_{k+1}\}}
\]
($X$ is supposed to be constant after time 1). Then $X^\varepsilon$ is
a martingale
in its filtration; if the proposition were proved for $X^\varepsilon$,
we would
have
%
%e3.17 ###
\begin{eqnarray}\label{vpeps}
V_p(X^\varepsilon)
&\le& C_p \mathbb{E}\sum_k|X_{\varepsilon
S_{k+1}}-X_{\varepsilon
S_k}
|^p\nonumber\\
&\le& C_p' \mathbb{E}\sum_k\Bigl(\sum|\Delta X_t|^21_{\{\varepsilon
S_k\le
t<\varepsilon
S_{k+1}\}}
\Bigr)^{p/2}\\
&\le& C_p' \mathbb{E}\sum|\Delta X_t|^p\nonumber
\end{eqnarray}
where we have used in the second line the classical Burkholder--Davis--Gundy
inequalities; it is then sufficient to let $\varepsilon$ tend to 0.
Thus it is
sufficient to prove the result for martingales varying only on a
sequence of
totally inaccessible stopping times. By separating the positive and negative
parts of the jumps, such a martingale is the difference of two martingales
with finite variation and with no negative jump, so we only have to prove
the result for these martingales.
\end{step}

\begin{step}\label{step2}
We suppose therefore that $X$ has finite
variation with positive jumps at a sequence of stopping times $S_k$. Thus
the positive part $dX^+=X^\nearrow$ of the Lebesgue--Stieltjes measure
of $X$
is purely atomic; it is carried by the times of jumps of $X$. Let $\tau
$ be
in $\mathbb{T}$; it is the projection of some $(\tau^\nearrow,x)$ of
$\mathbb{G}$,
and
\[
h(\tau)=\sup\{X_s-x; \tau^\nearrow\le s\le T(\tau^\nearrow
,x)\}
\]
with
\[
T(t,x):=\inf\{s\ge t;X_s\le x\}.
\]
Then \textup{(\ref{vphtau})} implies that
\begin{eqnarray*}
V_p(X) &\le& 2p\int_{\inf X}^{X_0}\sup\{X_s-x;s\le
T(0,x)\}^{p-1}\,dx
\\
&&{} +
2p\sum_{t\in J}\int_{X_{t-}}^{X_t}
\sup\{X_s-x;t\le s\le T(t,x)\}^{p-1}\,dx
\end{eqnarray*}
where $J=\{S_k; k\ge1\}$. The first term corresponds to the integral
on the
arc $[A,O]$ of $\mathbb{T}$, on which $\tau^\nearrow\le0$; its
expectation is
dominated by the expectation of $|X_1-X_0|^p$ (Doob's inequality) which can
be estimated by the right-hand side of \textup{(\ref{evp})} with the
technique of
\textup{(\ref{vpeps})}. The second term corresponds to the integral
on the remaining
part of the tree, for which $\tau^\nearrow\in J$. In order to
estimate it,
consider some jump $S=S_k$ and notice that since $X$ is a martingale
with no
negative jump,
\[
\mathbb{P}[\sup\{X_s-x;S\le s\le T(S,x)\}\ge a\mid\mathcal{F}
_S]
\le\frac{X_S-x}{a}
\]
for $X_{S-}\le x\le X_S$ and $a\ge X_S-x$. We deduce that
\[
\mathbb{E}[\sup\{X_s-x;S\le s\le T(S,x)\}
^{p-1}\mid\mathcal{F}
_S]
\le(X_S-x)\int_{X_S-x}^M a^{p-3}\,da,
\]
so
\[
\mathbb{E}\biggl[\int_{X_{S-}}^{X_S}
\sup\{X_s-x;S\le s\le T(t,x)\}^{p-1}\,dx\mid\mathcal{F}
_S\biggr]
\le(\Delta X_S)^p/\bigl(p(2-p)\bigr)
\]
and we can conclude by summing on the times of jumps $S=S_k$.\qed
\end{step}
\noqed\end{pf}

We now give for L\'{e}vy processes the analogue of Proposition \ref{dimfbm}.

\begin{proposition}
Let $X$ be an $\alpha$-stable L\'{e}vy process. Then, for almost any path
$\omega$ of $X$,
\[
\dim_H\mathbb{T}=\operatorname{\overline{dim}}\mathbb{T}=\mathcal
{V}(\omega)=\alpha\vee1.
\]
\end{proposition}

\begin{pf}
For $\alpha<1$, the process has finite variation, so $\mathcal
{V}(X)=1$ and
the dimension is 1. For $\alpha\ge1$, the fact that $\mathcal
{V}(X)\le
\alpha$
is classical and can be deduced from Proposition \ref{super}; thus
\[
1\le\dim_H\mathbb{T}\le\operatorname{\overline{dim}}\mathbb
{T}=\mathcal{V}(X)\le\alpha.
\]
Our result is therefore proved for $\alpha=1$. Suppose now $\alpha
>1$. We
will use the notation
\[
\delta(s,t)=\delta((s,X_s),(t,X_t))=X_s+X_t-2\inf_{[s,t]}X.
\]
It is known (Proposition VIII.2 of \cite{bertoin96}) that
\[
\mathbb{P}\biggl[\inf_{[0,1/2]}X>-u\biggr]\le C u^{\alpha\beta
}=O(u^{\alpha-1})
\]
as $u\downarrow0$, for $\beta=\mathbb{P}[X_t\le0]\ge(\alpha
-1)/\alpha
$. We also
have
\begin{eqnarray*}
\mathbb{P}[|X_1|<u\mid\mathcal{F}_{1/2}]
&\le&\sup_x\mathbb{P}[x-u<X_1-X_{1/2}<x+u]\\
&=&\sup_x\mathbb{P}[x-u<X_{1/2}<x+u]=O(u)
\end{eqnarray*}
because $X_{1/2}$ has a bounded density, so by taking the intersection of
these two events,
\[
\mathbb{P}[\delta(0,1)<u]=O(u^\alpha).
\]
We deduce that $\delta(0,1)^{-p}$ is integrable for any $p<\alpha$. The
variables $\delta(s,t)$ satisfy the same property, and by scaling,
\[
\mathbb{E}\delta(s,t)^{-p}=C(t-s)^{-p/\alpha}.
\]
This can be used to prove \textup{(\ref{energy})} for any $p<\alpha
$, so we
deduce as
in Proposition \ref{dimfbm} that the Hausdorff dimension is bounded
below by
$\alpha$.
\end{pf}

\begin{remark}
Another real tree, called the L\'{e}vy tree, has been associated to $X$ in
\cite{leglejan98a} when $X$ has only positive jumps. This tree is different
from $\mathbb{T}$ but is related to it; times which project on the same
point of
$\mathbb{T}$ also project on the same point of the L\'{e}vy tree, but an
arc of
$\mathbb{T}$ associated to a jump of $X$ is concentrated in the L\'{e}vy tree
into a
single point.
\end{remark}

Let us now give an analogue of Proposition \ref{nala} for L\'{e}vy processes.

\begin{proposition}\label{proposition314}
Let $X$ be a L\'{e}vy process. Suppose that almost surely, $X$~has no interval
on which it is monotone, and define
\[
\xi(a)=\mathbb{E}[S^a+T^a]
\]
for
\[
T^a:=\inf\biggl\{t;\>X_t<\sup_{[0,t]}X-a\biggr\},\qquad
S^a:=\inf\biggl\{t;\>X_t>\inf_{[0,t]}X+a\biggr\}.
\]
Then $\lim_0\xi=0$, and $\xi(a)N^a(X)$ (for the process $X$ on the time
interval $[0,1]$) converges in probability to 1 as $a\downarrow0$. If
$\xi(a)=O(a^\alpha)$ for some $\alpha>0$, then the convergence is almost
sure.
\end{proposition}

When the assumption about $X$ is not satisfied, then $X$ or $-X$ is the sum
of a subordinator and a compound Poisson process. In this case,
$\mathbb{T}
$ is
finite, so $N^a$ is bounded.

\begin{pf*}{Proof of Proposition \protect\ref{proposition314}}
Consider the times $T_i^a=T_i^a(X)$ and $S_i^a=S_i^a(X)$ defined by
\textup{(\ref{tiasia})}. On the other hand, notice that our
assumption implies that~$S^a$ and $T^a$ tend almost surely to 0 as $a\downarrow0$. Since
$X$ is
a L\'{e}vy process, times $T_{i+1}^a-S_i^a$ and $S_i^a-T_i^a$ are independent,
and have the same law as $T^a$ and~$S^a$. Thus
\[
\sup_{0\le t\le k\mu}\biggl(X_t-\inf_{[0,t]}X\biggr)\ge\sup_{1\le
j\le k}
\biggl(\sup_{(j-1)\mu\le t\le j\mu}
\biggl(X_t-\inf_{[(j-1)\mu,t]}X\biggr)\biggr)
\]
and the right-hand side is the supremum of $k$ independent identically
distributed variables, so
\[
\mathbb{P}[S^a>k\mu]
=\mathbb{P}\biggl[\sup_{0\le t\le
k\mu}\biggl(X_t-\inf_{[0,t]}X\biggr)<a\biggr]
\le(\mathbb{P}[S^a\ge\mu])^k
\]
for $\mu>0$. This probability is smaller than 1 from our assumption on $X$.
We deduce that the moments of $S^a$ (and $T^a$) are finite, so $\lim
_0\xi=0$
and
\[
\mathbb{P}[S^a>2k \mathbb{E}[S^a]]\le1/2^k.
\]
Thus $S^a/(2\mathbb{E}[S^a])$, and similarly $T^a/(2\mathbb
{E}[T^a])$, are
dominated by
a geometric variable, so the variances of $S^a$ and $T^a$ are dominated by
$(\mathbb{E}[S^a])^2$ and $(\mathbb{E}[T^a])^2$. Thus
%
%e3.18 ###
\begin{equation}\label{sna}\quad
\mathbb{E}[S_n^a]=n \xi(a),\qquad \operatorname{var}(S_n^a)=n
\bigl(\operatorname{var}(S^a)+\operatorname{var}(T^a)\bigr)\le C
n \xi(a)^2.
\end{equation}
If $n=n(a)\uparrow\infty$ as $a\downarrow0$, then
$n(a)^{-1}\xi(a)^{-1}S_{n(a)}$ has expectation 1 and has a variance
dominated by $1/n(a)$; in particular it converges in probability to 1. By
taking $n=n(a,\pm)\sim(1\pm\varepsilon)\xi(a)^{-1}$, we see from
\textup{(\ref{sna})} and
the definition of $N^a$ in \textup{(\ref{tiasia})} that $N^a$ is between
$n(a,-)$ and
$n(a,+)$ with a high probability, so the convergence in probability of the
proposition is proved. Moreover, for the second statement, it follows from the
Borel--Cantelli lemma that $n(a_k)^{-1}\xi(a_k)^{-1}S_{n(a_k)}$ converges
almost surely to 1 as soon as $\sum1/n(a_k)<\infty$. We can apply this
result to the above $n=n(a_k,\pm)$ for $a_k=1/k^\beta$ and $\beta$ large
enough, and we deduce that $\xi(a_k)N^{a_k}$ converges almost surely
to 1.
We conclude as in Proposition \ref{nala} from the monotonicity of $N^a$.
\end{pf*}

The almost sure convergence holds in particular for $\alpha$-stable
processes such that~$|X|$ is not a subordinator. In this case indeed,
$\xi(a)$ is proportional to $a^\alpha$ from the scaling property. For the
standard Brownian motion, $S^1$ and $T^1$ are the first hitting time of
1 by
a reflected Brownian motion, and have expectation 1. Thus $\xi
(a)=2a^2$ and
$N^a\sim1/(2a^2)$. This means that $C(1/2)=1/2$ in Proposition \ref{nala}.

We can deduce an estimation of $L^a$ when the process has infinite
variation. However, \textup{(\ref{lasum})} cannot be directly
applied; one has
to use
the associated continuous path, since times $S_i^a$ and $T_i^a$ can be jump
times.

%s4 ###
\section{Integrals and trees}
\label{integrals}

%s4.1 ###
\subsection{An integral on the tree}\label{sec41}

We now want to integrate some bounded function $\rho(t)$ against
$\omega$.
First suppose that $\omega$ is continuous and $\pi(0)=\pi(1)$. Let us
remember (Proposition \ref{finitevar}) that if $\omega$ has finite
variation, then $\omega^\nearrow$ and $\omega^\nwarrow$ are finite measures,
and are the positive and negative parts of the Lebesgue--Stieltjes measure
$d\omega$; moreover, since the images of the finite length measure
$\lambda$
by $\tau\mapsto\tau^\nearrow$ and $\tau\mapsto\tau^\nwarrow$ are
respectively $\omega^\nearrow$ and $\omega^\nwarrow$ (Proposition
\ref{lebesgue}), we have
\[
\int_0^1\rho(t)\omega^\nearrow(dt)=\int_\mathbb{T}\rho(\tau
^\nearrow
)\lambda
(d\tau),\qquad
 \int_0^1\rho(t)\omega^\nwarrow(dt)=\int_\mathbb{T}\rho
(\tau
^\nwarrow
)\lambda(d\tau),
\]
so
%
%e4.1 ###
\begin{equation}\label{int01}
\int_0^1\rho\,d\omega=\int_\mathbb{T}\bigl(\rho(\tau^\nearrow)
-\rho(\tau^\nwarrow)\bigr)\lambda(d\tau).
\end{equation}
If $\pi(0)\ne\pi(1)$, we extend $\omega$ to $[-1,2]$ as in
\textup{(\ref{extension})}, put $\rho(t)=0$ for $t\notin[0,1]$, and
we can
use the
same formula to define the integral; actually, with the notation
\textup{(\ref{newroot})}, we have
\begin{eqnarray}\label{intext}
\int_0^1\rho\,d\omega
&=&\int_{\mathbb{T}\setminus[A,B]}
\bigl(\rho
(\tau
^\nearrow)
-\rho(\tau^\nwarrow)\bigr)\lambda(d\tau)\nonumber\\[-8pt]\\[-8pt]
&&{}-\int_{[A,O]}\rho(\tau^\nwarrow)\lambda(d\tau)
+\int_{[O,B]}\rho(\tau^\nearrow)\lambda(d\tau).\nonumber
\end{eqnarray}
In this form, one can notice that the integral on $[0,1]$ depends on
$\rho$
and $\omega$ on $[0,1]$, and not on the extension of $\omega$ out of
$[0,1]$.

More generally, even if $\omega$ has infinite variation, we can define the
integral by the right-hand side of \textup{(\ref{int01})} or \textup
{(\ref{intext})}, provided
%
%e4.2 ###
\begin{equation}\label{integrable}
\int_\mathbb{T}|\rho(\tau^\nearrow)
-\rho(\tau^\nwarrow)|\lambda(d\tau)<\infty.
\end{equation}
Notice that the right-hand side of \textup{(\ref{int01})} is the
limit as
$a\downarrow0$ of the integral on the trimmed tree $\mathbb{T}^a$
which is the
tree of $\omega^a$ defined by \textup{(\ref{defomegaa})}, so $\int
\rho\,
d\omega$ is
the limit of $\int\rho\,d\omega^a$. This means that in this sense our
approach is similar to other approaches using a regularization of
$\omega$;
another example for which there has been a lot of work recently is the
Russo--Vallois approach \cite{rusval93}.

We now verify that we can apply our technique in the Young framework.

\begin{theorem}\label{young}
Assume that $\omega$ is continuous. One has
%
%e4.3 ###
\begin{equation}\label{inttabs}
\int_\mathbb{T}|\rho(\tau^\nearrow)-\rho(\tau^\nwarrow
)|\lambda
(d\tau)
\le C V_p(\omega)^{1/p}\bigl(V_q(\rho)^{1/q}+\sup|\rho|\bigr)
\end{equation}
for some $C=C(p,q)$, as soon as $1/p+1/q>1$. Thus \textup{(\ref
{integrable})} is
satisfied as soon as $1/\mathcal{V}(\omega)+1/\mathcal{V}(\rho)>1$, and
in this
case we can define $\int\rho\,d\omega$ by the right-hand side of
\textup{(\ref{int01})} or \textup{(\ref{intext})}. It satisfies
%
%e4.4 ###
\begin{equation}\label{estyoung}
\biggl|\int_0^1\rho\,d\omega\biggr|
\le C V_p(\omega)^{1/p}\bigl(V_q(\rho)^{1/q}+\sup|\rho|\bigr)
\end{equation}
for $1/p+1/q>1$. Moreover, this integral coincides with the
Riemann--Stieltjes integral constructed by Young \cite{young36} (see also
\cite{lyonsqian02,lyonscl07}); this means that
\[
\int_0^1\rho\,d\omega=\lim\sum_i\rho(s_i)\bigl(\omega
(t_{i+1})-\omega
(t_i)\bigr)
\]
for $t_i\le s_i\le t_{i+1}$, as the mesh of the subdivision $(t_i)$ of
$[0,1]$ tends to 0. The integral $\int_s^t\rho\,d\omega$ can be defined
similarly by replacing $\rho$ by $\rho 1_{(s,t]}$; it satisfies the Chasles
relation, and
%
%e4.5 ###
\begin{equation}\label{vpint}
V_p\biggl(\int_0^{\bolds{.}}\rho\,d\omega\biggr)^{1/p}
\le C V_p(\omega)^{1/p}\bigl(V_q(\rho)^{1/q}+\sup|\rho|\bigr).
\end{equation}
\end{theorem}

\begin{pf}
Let us first assume $\pi(0)=\pi(1)$. It follows from the disintegration
formula $\lambda=\lambda_2$ of Proposition \ref{lambda12} that
\[
I_a:=-\frac{\partial}{\partial a}\int_{\mathbb{T}^a}
|\rho(\tau^\nearrow)-\rho(\tau^\nwarrow)|\lambda
(d\tau)
=\sum_{\tau\in\partial\mathbb{T}^a}|\rho(\tau^\nearrow
)-\rho
(\tau
^\nwarrow)|.
\]
Define $0\le r<1$ by $1/q+r/p=1$. Then
%
%e4.6 ###
\begin{eqnarray}\label{iavpvq}
I_a&\le&\Biggl(\sum_{\tau\in\partial\mathbb{T}^a}
|\rho(\tau^\nearrow)-\rho(\tau^\nwarrow)|^q\Biggr)^{1/q}
(N^a)^{r/p}\nonumber\\[-8pt]\\[-8pt]
&\le&\frac1{2^{r/p}a^r}V_q(\rho)^{1/q}V_p(\omega)^{r/p}\nonumber
\end{eqnarray}
from H\"{o}lder's inequality and \textup{(\ref{vpapna})}. Consequently,
$I_a$ is of order
$1/a^r$ and is integrable with respect to $a$ near 0; more precisely, with
$\|\omega\|=\sup\omega-\inf\omega$,
%
%e4.7 ###
\begin{eqnarray}\label{estimint}\quad
\int_\mathbb{T}|\rho(\tau^\nearrow)-\rho(\tau^\nwarrow)
|\lambda(d\tau)
&=&\int_0^{\|\omega\|}I_a\,da\nonumber\\
&\le&\frac1{2^{r/p}(1-r)}\|\omega\|^{1-r}
V_q(\rho)^{1/q}V_p(\omega)^{r/p}\\
&\le&\frac1{2^{1/p}(1-r)}V_q(\rho)^{1/q}V_p(\omega)^{1/p}\nonumber
\end{eqnarray}
where we have used $V_p(\omega)\ge2\|\omega\|^p$ in the last line. If
$A=\pi(0)\ne\pi(1)=B$, we decompose $\mathbb{T}$ into $[A,B]$ and
$\mathbb{T}\setminus[A,B]$; we can apply the above procedure to the
integral on
the latter part, and again prove \textup{(\ref{estimint})}, but
without the
factor 2.
On the other hand, $[A,B]$ has finite length so the integral is finite on
it; more precisely,
\begin{eqnarray*}
\int_{[A,B]}|\rho(\tau^\nearrow)-\rho(\tau^\nwarrow)
|\lambda
(d\tau)
&=&\int_{[A,O]}|\rho(\tau^\nwarrow)|\lambda(d\tau)
+\int_{[O,B]}|\rho(\tau^\nearrow)|\lambda(d\tau)\\
&\le&\delta(0,1) \sup|\rho|\le2V_p(\omega)^{1/p} \sup|\rho|.
\end{eqnarray*}
The result \textup{(\ref{inttabs})} follows by adding these two
estimates. Thus we
can define the integral $\int_0^1\rho\,d\omega$ by \textup{(\ref
{int01})}; this
integral satisfies \textup{(\ref{estyoung})}, and similarly,
\[
\biggl|\int_s^t\rho\,d\omega\biggr|
\le C V_p(\omega;s,t)^{1/p}\bigl(V_q(\rho)^{1/q}+\sup|\rho|\bigr)
\]
where the $p$-variation of $\omega$ is limited to $[s,t]$. One easily
deduces \textup{(\ref{vpint})} by applying
%
%e4.8 ###
\begin{equation}\label{sumivp}
\sum_iV_p(\omega;t_i,t_{i+1})\le V_p(\omega).
\end{equation}
More precisely, by considering the variations of $\rho$ and $\omega$ on
$[s,t]$,
\begin{eqnarray*}
&&\biggl|\int_s^t\rho\,d\omega-\rho(s)\bigl(\omega(t)-\omega
(s)\bigr)\biggr|\\
&&\qquad=\biggl|\int_s^t\bigl(\rho(\cdot)-\rho(s)\bigr)\,d\omega\biggr|\\
&&\qquad\le C V_p(\omega;s,t)^{1/p}\bigl(V_q(\rho;s,t)^{1/q}
+\sup|\rho(\cdot)-\rho(s)|\bigr)\\
&&\qquad\le C' V_p(\omega;s,t)^{1/p}V_q(\rho;s,t)^{1/q}.
\end{eqnarray*}
Thus
%
%e4.9 ###
\begin{eqnarray}\label{cprime}
&&\biggl|\int_{t_i}^{t_{i+1}}\rho\,d\omega-\rho(s_i)
\bigl(\omega(t_{i+1})-\omega(t_i)\bigr)\biggr|\nonumber\\
&&\qquad\le C V_q(\rho;t_i,t_{i+1})^{1/q}V_p(\omega
;t_i,t_{i+1})^{1/p}
\\
&&\qquad\le C'\bigl(V_q(\rho;t_i,t_{i+1})+V_p(\omega;t_i,t_{i+1})\bigr)
V_p(\omega;t_i,t_{i+1})^{(1-r)/p}.\nonumber
\end{eqnarray}
By applying \textup{(\ref{sumivp})} and the similar estimate for
$\rho$, we get
\begin{eqnarray*}
&&\Biggl|\int_0^1\rho\,d\omega-\sum_i\rho(s_i)
\bigl(\omega(t_{i+1})-\omega(t_i)\bigr)\Biggr|\\
&&\qquad\le C\bigl(V_q(\rho)+V_p(\omega)\bigr)\sup_iV_p(\omega
;t_i,t_{i+1})^{(1-r)/p}
\end{eqnarray*}
which converges to 0 since $\omega$ is continuous.
\end{pf}

\begin{remark}
In the proof, we have considered separately the arc $[A,B]$. Actually,
%
%e4.10 ###
\begin{equation}\label{intab}
\int_{[A,B]}\bigl(\rho(\tau^\nearrow)-\rho(\tau^\nwarrow)
\bigr)\lambda
(d\tau)
=\int\rho\,d\underline{\omega}
\end{equation}
with
%
%e4.11 ###
\begin{equation}\label{undomega}
\underline{\omega}(t)=\inf_{[0,t]}\omega\vee\inf_{[t,1]}\omega.
\end{equation}
\end{remark}

\begin{remark}
In the framework of Theorem \ref{young}, the fact that our integral is a
Riemann--Stieltjes integral implies that it is linear with respect to
$\omega$; this property was not evident on our definition, since the tree
associated to the sum of two paths is not simply related to the trees
of the
two paths. Actually, we do not know whether the space of $\omega$ satisfying
\textup{(\ref{integrable})} is linear.
\end{remark}

\begin{remark}
Young integrals can also be written as classical integrals on the time
interval by means of a completely different technique, namely fractional
differential calculus (see \cite{zahle98}).
\end{remark}

\begin{theorem}\label{cadlag}
Theorem \ref{young} holds for c\`{a}dl\`{a}g paths $\omega$, provided $\rho$ is
continuous at times of discontinuity of $\omega$.
\end{theorem}

\begin{pf}
The tree is associated to a continuous path
$(\omega'(t,x); (t,x)\in\mathbb{G})$, as it has been explained in
Proposition
\ref{graphtree}, and $\omega'$ has the same variations as $\omega
$. One
can also consider $\rho'(t,x)=\rho(t)$ which has the same variations as
$\rho$. Then the left-hand side of \textup{(\ref{inttabs})} is the
integral for
$\rho'$ and $\omega'$, so \textup{(\ref{inttabs})} holds true. For the
Riemann sums,
we modify \textup{(\ref{cprime})} in the previous proof by introducing
$r'<1$ such that
$1/p+1/q=1/r'$; then
\begin{eqnarray*}
&&\biggl|\int_{t_i}^{t_{i+1}}\rho\,d\omega-\rho(s_i)
\bigl(\omega(t_{i+1})-\omega(t_i)\bigr)\biggr|\\
&&\qquad\le
C\bigl(V_q(\rho;t_i,t_{i+1})+V_p(\omega;t_i,t_{i+1})\bigr)\\
&&\qquad\quad{}\times
V_p(\omega;t_i,t_{i+1})^{(1-r')/p}V_q(\rho;t_i,t_{i+1})^{(1-r')/q},
\end{eqnarray*}
so that
\begin{eqnarray*}
&&\Biggl|\int_0^1\rho\,d\omega-\sum_i\rho(s_i)
\bigl(\omega(t_{i+1})-\omega(t_i)\bigr)\Biggr|\\
&&\qquad \le C\bigl(V_q(\rho)+V_p(\omega)\bigr)
\sup_i(V_p(\omega;t_i,t_{i+1})^{1/p}V_q(\rho
;t_i,t_{i+1})^{1/q})^{1-r'}.
\end{eqnarray*}
We have to prove that the supremum tends to 0 as the mesh of the subdivision
tends to 0. For any $\varepsilon>0$, let us consider
\[
J_\varepsilon:=\{i; |\Delta\omega(t)|\ge
\varepsilon\mbox{
for some }
t_i<t\le t_{i+1}\}.
\]
Then
\[
\limsup\>\sup_{i\notin J_\varepsilon}V_p(\omega
;t_i,t_{i+1})^{1/p}\le\varepsilon,
\]
and the number of jumps greater than $\varepsilon$ is finite, so from the
continuity of $\rho$ at these points,
\[
\lim\>\sup_{i\in J_\varepsilon}V_q(\rho;t_i,t_{i+1})^{1/q}=0.
\]
We deduce the convergence from these two properties.
\end{pf}

\begin{remark}
If $\rho$ and $\omega$ have common discontinuity times, our integral can
still be defined, but the Riemann--Stieltjes approach has to be
modified, as
in the classical Young work \cite{young36}.
\end{remark}

This theory can be applied to paths of fractional Brownian motions with
Hurst parameter $H>1/2$, or to L\'{e}vy processes without Brownian part and such
that $|x|^p\wedge1$ is integrable with respect to the L\'{e}vy measure for some
$p<2$.

%s4.2 ###
\subsection{Beyond the Young integral}\label{sec42}

A limitation of the Young integral concerns its iteration. If $\omega$ and
$\rho$ have respectively $p$- and $q$-finite variation for $1/p+1/q>1$ and
$\omega$ is continuous, then we can consider the function
\[
x(t):=x_0+\int_0^t\rho\,d\omega,
\]
and \textup{(\ref{vpint})} implies that $x$ has $p$-finite variation;
however, it
generally does not have $q$-finite variation so, unless $p<2$, one cannot
construct $\int x\,d\omega$. Nevertheless, we now check that this is
possible with our framework (for a continuous one-dimensional path
$\omega$). The idea is to look for a weaker condition than
$V_q(\rho)<\infty$ for \textup{(\ref{integrable})}.

For instance, if $\rho(t)=f(\omega(t))$, \textup{(\ref
{integrable})} holds
for any
bounded $f$ and any continuous $\omega$, and
\[
\int_0^1f(\omega(t))\,d\omega(t)=F(\omega(1))-F(\omega(0))
\]
for a primitive function $F$ of $f$; this is because the integral on
$\mathbb{T}\setminus[A,B]$ in~\textup{(\ref{intext})} is 0
($f(\omega(\tau^\nearrow))=f(\omega(\tau^\nwarrow))$), and the
integral on
$[A,B]$ is easily computed from~\textup{(\ref{intab})}. However, in
this case, the
integral is not always the limit of Riemann sums, as it is easily seen for
$f(x)=x$. We want to generalize this example.

Define
\[
V_q(\rho|\omega):=\sup\sum_k|\rho(t_{2k+2})-\rho
(t_{2k+1})|^q
\]
where the supremum is with respect to subdivisions $(t_i)$ of $[0,1]$ such
that $\omega(t_{2k+1})=\omega(t_{2k+2})$, and put
\[
\mathcal{V}(\rho|\omega):=\inf\{q\ge1; V_q(\rho|\omega
)<\infty
\}
\le\mathcal{V}(\rho).
\]

\begin{theorem}\label{eyoung}
Let $\omega$ be continuous. The integrability condition \textup{(\ref
{integrable})}
holds as soon as
\[
1/\mathcal{V}(\omega)+1/\mathcal{V}(\rho|\omega)>1.
\]
Moreover, if $1/p+1/q>1$ and $\omega$ fixed with $V_p(\omega)<\infty
$, the
space of bounded functions $\rho$ such that $V_q(\rho|\omega)<\infty$
is a
Banach space $\mathbb{B}_{q,\omega}$ for the norm
\[
\|\rho\|_{q,\omega}:=V_q(\rho|\omega)^{1/q}+\sup|\rho|,
\]
and we have
%
%e4.12 ###
\begin{eqnarray}\label{vqint}
V_q\biggl(\int_0^{\bolds{.}}\rho\,d\omega\bigm|\omega\biggr)^{1/q}
&\le& C V_p(\omega)^{1/p}V_q(\rho|\omega)^{1/q},
\\
\label{normint}
\biggl\|\int_0^{\bolds{.}}\rho\,d\omega\biggr\|_{q,\omega}+V_p\biggl(\int
_0^{\bolds{.}}\rho\,
d\omega\biggr)^{1/p}
&\le& C V_p(\omega)^{1/p}\|\rho\|_{q,\omega},
\end{eqnarray}
for some $C=C(p,q)$.
\end{theorem}

\begin{pf}
In the\vspace*{2pt} estimation \textup{(\ref{iavpvq})}, we can use $V_q(\rho
|\omega)$
instead of
$V_q(\rho)$ since we consider the subdivisions defined by
$t_{2k+1}=\tau^\nearrow$ and $t_{2k+2}=\tau^\nwarrow$ for
$\tau\in\partial\mathbb{T}^a$. Thus \textup{(\ref{inttabs})} is
replaced by
%
%e4.14 ###
\begin{equation}\label{estimvpomrho}
\int_\mathbb{T}|\rho(\tau^\nearrow)-\rho(\tau^\nwarrow
)|\lambda
(d\tau)\le
C V_p(\omega)^{1/p}\|\rho\|_{q,\omega}.
\end{equation}
This proves the first statement. The Banach property is easily verified from
the lower semicontinuity of $\rho\mapsto V_q(\rho|\omega)$ with
respect to
uniform convergence. By applying \textup{(\ref{estimvpomrho})} on
$[s,t]$, we
estimate $\int_s^t\rho\,d\omega$, and deduce that $\int_0^{\bolds{.}}\rho\,
d\omega$
and $V_p(\int\rho\,d\omega)^{1/p}$ are bounded by the right-hand
side of
\textup{(\ref{estimvpomrho})} [for the estimation of the
$p$-variation, we use
\textup{(\ref{sumivp})}]. The last property which has to be proved in
order to
conclude is \textup{(\ref{vqint})}. To this end, we are going to
check that
%
%e4.15 ###
\begin{equation}\label{intbar}
\biggl|\int_0^1\rho\,d\omega\biggr|\le C V_p(\omega)^{1/p}
V_q(\rho|\omega)^{1/q}
\end{equation}
as soon as $\omega(0)=\omega(1)$; then \textup{(\ref{vqint})}
follows by applying
\textup{(\ref{intbar})} on the intervals $[t_{2k+1},t_{2k+2}]$ in
order to estimate
$V_q(\cdot|\omega)$. The left-hand side of \textup{(\ref{intbar})} is
written as an
integral on the tree; the integral on $\mathbb{T}\setminus[A,B]$ is
estimated by
the right-hand side of \textup{(\ref{intbar})} as in \textup{(\ref
{estimint})}; for the
integral on $[A,B]$, it can be written as
\[
\int_{[A,B]}\bigl(\rho(\tau^\nearrow)-\rho(\tau^\nwarrow)
\bigr)\lambda
(d\tau)
=\int_{\inf\omega}^{\omega(0)}\bigl(\rho(\beta_2(x))-\rho(\beta
_1(x))\bigr)\,dx
\]
with
\[
\beta_1(x)=\inf\{t;\omega(t)=x\},\qquad
\beta_2(x)=\sup\{t;\omega(t)=x\}.
\]
This expression is also easily estimated by the right-hand side of
\textup{(\ref{intbar})}.
\end{pf}

As an application, we can solve differential equations driven by a
multidimensional path, provided all the components of the path but one are
smooth enough.

\begin{theorem}\label{edo}
For $1/p+1/q>1$ and $q\le p$, consider a continuous real-valued map
$\omega$
with finite $p$-variation, and let $\mathbb{B}_{p,q,\omega}$ be the Banach
space of functions $\rho$ such that
\[
\|\rho\|_{p,q,\omega}:=V_p(\rho)^{1/p}+V_q(\rho|\omega)^{1/q}+\sup
|\rho|
\]
is finite. Consider also a continuous function $\eta$ with values in
$\mathbb{R}^{d-1}$ and with finite $q$-variation, and let $\xi
=(\omega
,\eta)$
with values in $\mathbb{R}^d$. Let $f$ be a $C^2$ function with bounded
derivatives from $\mathbb{R}^n$ into the space of linear maps
$\mathcal{L}(\mathbb{R}^d,\mathbb{R}^n)$. Consider, for $x_0$ in
$\mathbb{R}^n$, the equation
\[
x(t)=x_0+\int_0^tf(x(s))\,d\xi(s)
\]
where the integral should be understood as the sum of integrals with respect
to each component, each one being given by an expression of type
\textup{(\ref{int01})} or \textup{(\ref{intext})}. Then this
equation has a unique
solution in
the Banach space $\mathbb{B}_{p,q,\omega}^n$.
\end{theorem}

\begin{pf}
In this proof, the constants $C$ may depend on $f$ and $x_0$, but not on
$\xi$. It is not difficult to deduce from the Lipschitz property of
$f$ that
\[
F\dvtx\bigl(x(t); 0\le t\le1\bigr)\mapsto
\bigl(f(x(t)); 0\le t\le1\bigr)
\]
maps $(\mathbb{B}_{p,q,\omega})^n$ into $(\mathbb{B}_{p,q,\omega})^{nd}$
and has
at most linear growth:
%
%e4.16 ###
\begin{equation}\label{fxq}
\|F(x)\|_{p,q,\omega}\le C(\|x\|_{p,q,\omega
}+1).
\end{equation}
Let us prove that $F$ is locally Lipschitz. It is easy to verify
%
%e4.17 ###
\begin{equation}\label{supff}
\sup|F(x_2)-F(x_1)|\le C \sup|x_2-x_1|,
\end{equation}
and let us estimate $V_q(F(x_2)-F(x_1)|\omega)$. Let $(t_i)$ be a
subdivision satisfying $\omega(t_{2k+1})=\omega(t_{2k+2})$, and use the
notation $\Delta_iv=v(t_{i+1})-v(t_i)$. It follows from the
boundedness of
the derivatives of $f$ that
\begin{eqnarray*}
&&|f(x_2(t_{i+1}))-f(x_1(t_{i+1}))
-f(x_2(t_i))+f(x_1(t_i))|\\
&&\qquad\le C\bigl(|x_2(t_{i+1})-x_1(t_{i+1})|+
|x_2(t_i)-x_1(t_i)|\bigr)\\
&&\qquad\quad{}\times
(|\Delta_ix_2|+|\Delta_ix_1|)
+ C |\Delta_ix_2-\Delta_ix_1|\\
&&\qquad \le2C \sup|x_2-x_1|
(|\Delta_ix_2|+|\Delta_ix_1|)+C
|\Delta_ix_2-\Delta_ix_1|.
\end{eqnarray*}
By taking the $q$th power and summing over indices $i=2k+1$, we deduce
%
%e4.18 ###
\begin{eqnarray}\label{vqff}
&& V_q\bigl(F(x_2)-F(x_1)\mid\omega\bigr)\nonumber\\
&&\qquad \le C \sup|x_2-x_1|^q
\bigl(V_q(x_1|\omega)+V_q(x_2|\omega)\bigr)
+C V_q(x_2-x_1|\omega)\\
&&\qquad \le C \|x_2-x_1\|_{p,q,\omega}^q
\bigl(V_q(x_1|\omega)+V_q(x_2|\omega)+1\bigr).\nonumber
\end{eqnarray}
We prove similarly that
%
%e4.19 ###
\begin{equation}\label{vpff}
V_p\bigl(F(x_2)-F(x_1)\bigr)\le C \|x_2-x_1\|
_{p,q,\omega}^p
\bigl(V_p(x_1)+V_p(x_2)+1\bigr).
\end{equation}
It follows from \textup{(\ref{supff})}, \textup{(\ref{vqff})} and
\textup{(\ref{vpff})} that
$F$ is
locally Lipschitz; more precisely,
%
%e4.20 ###
\begin{equation}\label{fx2x1}\qquad\quad
\|F(x_2)-F(x_1)\|_{p,q,\omega}\le C(\|x_1\|
_{p,q,\omega}
+\|x_2\|_{p,q,\omega}+1)\|x_2-x_1\|_{p,q,\omega}.
\end{equation}
On the other hand, the property $q\le p$ and Theorem \ref{young} (applied
with an exchange of $p$ and $q$) show that
\begin{eqnarray*}
\biggl\|\int_0^{\bolds{.}}\rho\,d\eta\biggr\|_{p,q,\omega}
\le C \sup\biggl|\int_0^{\bolds{.}}\rho\,d\eta\biggr|+C V_q\biggl(\int
_0^{\bolds{.}}\rho\,
d\eta\biggr)^{1/q}
\le C'\|\rho\|_{p,q,\omega}V_q(\eta)^{1/q}
\end{eqnarray*}
if $\rho$ takes its values in $\mathcal{L}(\mathbb{R}^{d-1},\mathbb
{R}^n)$. If now
$\rho$
takes its values in $\mathcal{L}(\mathbb{R}^d,\mathbb{R}^n)$, we
deduce by using also
\textup{(\ref{normint})} that
%
%e4.21 ###
\begin{equation}\label{normvec}
\biggl\|\int_0^{\bolds{.}}\rho\,d\xi\biggr\|_{p,q,\omega}
\le C\|\rho\|_{p,q,\omega}\bigl(V_p(\omega)^{1/p}+V_q(\eta
)^{1/q}\bigr).
\end{equation}
Thus, by joining \textup{(\ref{fxq})}, \textup{(\ref{fx2x1})} and
\textup{(\ref{normvec})}, we obtain
that the map
\[
\Phi\dvtx\bigl(\rho(t); 0\le t\le1\bigr)\mapsto
\biggl(x_0+\int_0^t\rho\,d\xi; 0\le t\le1\biggr)
\]
satisfies
\[
\|(\Phi\circ F)(x)\|_{p,q,\omega}\le
C+C\bigl(V_p(\omega)^{1/p}+V_q(\eta)^{1/q}\bigr)(1+\|x\|
_{p,q,\omega})
\]
and
\begin{eqnarray*}
&&\|(\Phi\circ F)(x_2)-(\Phi\circ F)(x_1)\|_{p,q,\omega}\\
&&\qquad\le C\bigl(V_p(\omega)^{1/p}+V_q(\eta)^{1/q}\bigr)(\|x_1\|
_{p,q,\omega}
+\|x_2\|_{p,q,\omega}+1)\|x_2-x_1\|_{p,q,\omega}.
\end{eqnarray*}
It is then classical to deduce that $\Phi\circ F$ has a unique fixed point
if $V_p(\omega)$ and $V_q(\eta)$ are small enough. We conclude like for
usual differential equations by dividing $[0,1]$ into subintervals where
$\omega$ and $\eta$ have small variation.
\end{pf}

In particular, we can work out a calculus for one-dimensional fractional
Brownian motions of any Hurst parameter, and the stochastic integrals
can be
interpreted as integrals on the tree; another interpretation can be worked
out by modifying Russo--Vallois integrals \cite{gradnouruva05,nourdin07}.

%s4.3 ###
\subsection{Integration for fractional Brownian motion}\label{sec43}

Up to now, we have found sufficient conditions ensuring that the integral
$\int\rho\,d\omega$ can be defined as an integral on the tree.
However, by
means of the disintegration $\lambda=\lambda_2$ of the length measure
(Proposition \ref{lambda12}), the strong integrability condition
\textup{(\ref{integrable})} can be replaced by the weaker condition
%
%e4.22 ###
\begin{equation}\label{wintegrable}
\int\Biggl|\sum_{\tau\in\partial\mathbb{T}^a}\bigl(
\rho(\tau^\nearrow)-\rho(\tau^\nwarrow)\bigr)\Biggr|\,da<\infty
\end{equation}
(where the number of terms in the sum is finite), and in this case we can
define
%
%e4.23 ###
\begin{equation}\label{wint01}
\int_0^1\rho\,d\omega:=\int\sum_{\tau\in\partial\mathbb
{T}^a}\bigl(
\rho(\tau^\nearrow)-\rho(\tau^\nwarrow)\bigr)\,da
\end{equation}
[with a form similar to \textup{(\ref{intext})} if $\pi(0)\ne\pi
(1)$]. This
is a
generalization of the previous framework, and the integral, when it exists, is
again the limit of $\int\rho\,d\omega^a$. If \textup{(\ref
{wintegrable})} is
satisfied for $\rho$ replaced by 0 out of $[s,t]$, we can define similarly
$\int_s^t\rho\,d\omega$ satisfying the Chasles relation. Our aim is
now to
check that this integral is well adapted to the differential calculus with
respect to a finite-dimensional $H$-fractional Brownian motion, for
$1/3<H\le1/2$ (made of independent one-dimensional fractional Brownian
motions), and that the integrals coincide with those of the rough paths
theory \cite{lyons98,lyonsqian02,lyonscl07,couqian02,lejay03}. Some
related results for the standard Brownian case $H=1/2$ are also given in
\cite{picard06}; in this case, the integrals which are considered here are
Stratonovich integrals, but it is also explained in \cite{picard06}
how one
can use the tree $\mathbb{T}$ to obtain It\^{o} integrals. We are going to consider
the two-dimensional case (higher dimension is similar).

\begin{theorem}\label{intfbm}
Consider a two-dimensional $H$-fractional Brownian motion for $H\le1/2$.
Then almost any path $(\omega,\eta)$ satisfies the following
properties:\vspace*{-5pt}
\begin{longlist}[2.]
\item[1.] Suppose $H>1/4$ and let $1/4<r<H$. Then the integral
$\int_s^t\eta\,d\omega$ can be defined in the sense of
\textup{(\ref{wint01})}. Moreover
%
%e4.24 ###
\begin{equation}\label{gammast}
\gamma(s,t):=\int_s^t\eta\,d\omega-\eta(s)\bigl(\omega(t)-\omega
(s)\bigr)
\end{equation}
satisfies
\[
|\gamma(s,t)|\le K(t-s)^{2r},
\]
where $K$ depends on $r$ and the path $(\omega,\eta)$, but not on $(s,t)$.
\item[2.]
Suppose $H>1/3$ and let $1/3<r<H$. Let $\rho$, $\phi$ and $\psi$ be bounded
paths such that
%
%e4.25 ###
\begin{equation}\label{hyprho}\qquad\quad
\bigl|\rho(t)-\rho(s)-\phi(s)\bigl(\omega(t)-\omega(s)\bigr)
-\psi(s)\bigl(\eta(t)-\eta(s)\bigr)\bigr|
\le K_1(t-s)^{2r}
\end{equation}
and
\[
|\psi(t)-\psi(s)|\le K_2(t-s)^r
\]
for any $s<t$ [where $K_1$ and $K_2$ may depend on $(\omega,\eta)$].
Then the integral $\int_s^t\rho\,d\omega$ can be defined in the
sense of
\textup{(\ref{wint01})}, and
\begin{eqnarray}\label{expansion}\qquad
&&\biggl|\int_s^t\rho\,d\omega-\rho(s)\bigl(\omega(t)-\omega
(s)\bigr)
-\frac{\phi(s)}2\bigl(\omega(t)-\omega(s)\bigr)^2-\psi(s)\gamma
(s,t)\biggr|\nonumber\\[-8pt]\\[-8pt]
&&\qquad \le K_3(t-s)^{3r}.\nonumber
\end{eqnarray}
\end{longlist}
\end{theorem}

\begin{pf}
Let $\mathbb{E}^\omega$ denote the integration with respect to the
law of
$\eta$,
with $\omega$ fixed, and let $\mathbb{T}$ be the tree of $\omega$. We
divide the
proof of the two parts of the theorem into two steps.

\begin{stepp}\label{stepp1}
Define a process $U^a$ as follows: consider the
points $\tau_1$, $\tau_2,\ldots$ of $\partial\mathbb{T}^a$ such that
$[\tau_i^\nearrow,\tau_i^\nwarrow]\subset[0,1]$, and let $U^a$ be
0 before
$\tau_1^\nearrow$, be constant on each
$[\tau_i^\nwarrow,\tau_{i+1}^\nearrow]$ and after the last
$\tau_i^\nwarrow$, be affine on each $[\tau_i^\nearrow,\tau
_i^\nwarrow]$,
and have the same increment on this interval as $\eta$. We will use the
notation $\Delta\tau_i=\tau_i^\nwarrow-\tau_i^\nearrow$. Since the
increments of $\eta$ are negatively correlated, we have for $j\le k$ and
$\varepsilon$ small enough
\begin{eqnarray*}
\mathbb{E}^\omega\bigl(U^a(\tau_k^\nwarrow)-U^a(\tau_j^\nearrow
)\bigr)^2
&\le& C\sum_{i=j}^k(\Delta\tau_i)^{2H}\\
&\le& C\biggl(\inf_i\Delta\tau_i\biggr)^{-2H+2\varepsilon}
\sum_{i=j}^k(\Delta\tau_i)^{4H-2\varepsilon}\\
&\le& C\biggl(\inf_i\Delta\tau_i\biggr)^{-2H+2\varepsilon}
\Biggl(\sum_{i=j}^k\Delta\tau_i\Biggr)^{4H-2\varepsilon}\\
&\le& K a^{-2+2\varepsilon}(\tau_k^\nwarrow-\tau_j^\nearrow
)^{4H-2\varepsilon},
\end{eqnarray*}
with $K=K(\omega)$ bounded in the spaces $L^q$; in the last line, we have
used the modulus of continuity of $\omega$. Thus $U^a$ is H\"{o}lder continuous
in $L^2(\mathbb{P}^\omega)$ on the set of times
$\{\tau_j^\nearrow,\tau_j^\nwarrow\}$; since it is extended to
$[0,1]$ by
affine interpolation, it satisfies the same property on the whole interval,
so
\[
\mathbb{E}^\omega\bigl(U^a(t)-U^a(s)\bigr)^2
\le K a^{-2+2\varepsilon}(t-s)^{4H-2\varepsilon}.
\]
Since the variable is conditionally Gaussian, estimates in
$L^q(\mathbb{P}^\omega)$ can be deduced for any $q$, so that, after
integration
with respect to $\omega$,
\[
\|U^a(t)-U^a(s)\|_{L^q}\le C_qa^{-1+\varepsilon
}(t-s)^{2H-\varepsilon}.
\]
By applying the Kolmogorov lemma,
%
%e4.26 ###
\begin{equation}\label{uaka}
|U^a(t)-U^a(s)|\le K^a a^{-1+\varepsilon}(t-s)^{2r}
\end{equation}
with $K^a$ bounded in $L^q$, uniformly in $a$, and for
$1/4<r<H-\varepsilon/2$.
Moreover,
\[
I_{s,t}^a:=\sum_{\tau\in\partial\mathbb{T}^a\dvtx[\tau^\nearrow,\tau
^\nwarrow
]\subset[s,t]}
\bigl(\eta(\tau^\nwarrow)-\eta(\tau^\nearrow)\bigr)
\]
is an increment of $U^a$ on a subinterval of $[s,t]$, so
\[
|I_{s,t}^a|\le K^a a^{-1+\varepsilon}(t-s)^{2r}.
\]
Since $K^a$ is bounded in $L^1$, $\int_0^{a_0}K^a a^{-1+\varepsilon
}\,da$ is finite
for any $a_0$ and almost any $(\omega,\eta)$; moreover, $I_{s,t}^a$ is
0 if
$a$ is greater than the oscillation of $\omega$. Thus
%
%e4.27 ###
\begin{equation}\label{intista}
\int_0^\infty|I_{s,t}^a|\,da\le K (t-s)^{2r}
\end{equation}
for some finite variable $K$. This implies that the integral
$\int_s^t\eta\,d\omega$ is well defined as claimed in the theorem, and
\[
\int_s^t\eta\,d\omega=\int_0^\infty I_{s,t}^a\,da
+\int_s^t\eta\,d\underline{\omega}
\]
where $[0,1]$ is replaced by $[s,t]$ in the notation \textup{(\ref
{undomega})}. The
estimation of $\gamma(s,t)$ follows from \textup{(\ref{intista})}
and the
moduli of
continuity of $\omega$ and $\eta$.
\end{stepp}

\begin{stepp}\label{stepp2}
Let us now consider the integral of $\rho$.
As in the
previous step, we consider the term $K^a$ of \textup{(\ref{uaka})},
and a path
$(\omega,\eta)$ such that $\int_0^{a_0}K^aa^{-1+\varepsilon}\,da$ is
finite. Consider
as in the previous step the times $\tau_j$ of $\partial\mathbb{T}^a$,
and define
$\psi^a$ to be $\psi_{\tau_j^\nearrow}$ on each
$[\tau_j^\nearrow,\tau_{j+1}^\nearrow)$, and 0 before $\tau
_1^\nearrow$;
then, by limiting the sums to indices $j$ such that
$[\tau_j^\nearrow,\tau_j^\nwarrow]\subset[s,t]$,
%
%e4.28 ###
\begin{eqnarray}\label{jsta}\qquad
J_{s,t}^a:\!&=&\sum_j
\bigl(\rho(\tau_j^\nwarrow)-\rho(\tau_j^\nearrow)\bigr)\nonumber\\[-8pt]\\[-8pt]
&=&\int_{s'}^{t'}\psi^a\,dU^a+\sum_j
\bigl(\rho(\tau_j^\nwarrow)-\rho(\tau_j^\nearrow)-\psi(\tau
_j^\nearrow)
\bigl(\eta(\tau_j^\nwarrow)-\eta(\tau_j^\nearrow)\bigr)\bigr)\nonumber
\end{eqnarray}
where $s'$ and $t'$ are the first $\tau_j^\nearrow$ and the last
$\tau_j^\nwarrow$ in $[s,t]$. Since $1/r+1/(2r)>1$, the first term is
estimated as a Young integral by means of \textup{(\ref{estyoung})}, so
%
%e4.29 ###
\begin{eqnarray}\label{intsprime}
&&\biggl|\int_{s'}^{t'}\psi^a\,dU^a\biggr|
\le C V_{1/(2r)}(U^a)^{2r}\bigl(V_{1/r}(\psi^a)^r+\sup|\psi
^a|\bigr)\nonumber\\[-8pt]\\[-8pt]
&&\qquad\le K K^aa^{-1+\varepsilon}(t-s)^{2r}\nonumber
\end{eqnarray}
for a finite $K$, and for $K^a$ obtained in the previous step. The second term
of \textup{(\ref{jsta})} is dominated from \textup{(\ref{hyprho})} by
%
%e4.30 ###
\begin{equation}\label{sumtauj}\quad
\sum(\tau_j^\nwarrow-\tau_j^\nearrow)^{2r}
\le K a^{-1+\varepsilon}\sum(\tau_j^\nwarrow-\tau_j^\nearrow)^{3r}
\le K a^{-1+\varepsilon}(t-s)^{3r}
\end{equation}
where we have used the modulus of continuity of $\omega$ in the first
inequality. Thus, by adding \textup{(\ref{intsprime})} and \textup
{(\ref{sumtauj})}, the
expression $J_{s,t}^a$ of \textup{(\ref{jsta})} is integrable with
respect to~$a$,
and $\int_s^t\rho\,d\omega$ is defined. Moreover,
\[
\int_s^t\rho\,d\omega=\int_s^t\rho\,d\underline{\omega}+\int J_{s,t}^a\,da.
\]
If $\rho(s)=\phi(s)=\psi(s)=0$, then $\rho$ is at most of order
$(t-s)^{2r}$, so the first term is at most of order $(t-s)^{3r}$; on the
other hand, in this case, one can put the exponent~$3r$ instead of $2r$ in
\textup{(\ref{intsprime})}, so the integral of $J_{s,t}^a$ is also of order
$(t-s)^{3r}$; thus $\int_s^t\rho\,d\omega$ is of order $(t-s)^{3r}$. This
can be applied to the integral of
\[
\rho(\cdot)-\rho(s)-\phi(s)\bigl(\omega(\cdot)-\omega(s)\bigr)-\psi(s)\bigl(\eta
(\cdot)-\eta(s)\bigr),
\]
and we deduce \textup{(\ref{expansion})}.\qed
\end{stepp}
\noqed\end{pf}

\begin{remark}
The estimate \textup{(\ref{expansion})} shows that the integral can be
constructed by
time discretization as limits of generalized Riemann sums
\begin{eqnarray}\label{taylor}\qquad
\int_0^1\rho\,d\omega&=&\lim\sum_i\biggl(\rho(t_i)\bigl(\omega
(t_{i+1})-\omega(t_i)\bigr)\nonumber\\[-8pt]\\[-8pt]
&&\hspace*{35.1pt}{}
+\frac{\phi(t_i)}{2}\bigl(\omega(t_{i+1})-\omega(t_i)
\bigr)^2+\psi(t_i)\gamma(t_i,t_{i+1})\biggr).\nonumber
\end{eqnarray}
\end{remark}

In the framework of Theorem \ref{intfbm}, we can construct similarly
integrals with respect to $\eta$ by means of the tree of $\eta$. Let
$(e_1,e_2)$ be the canonical basis of $\mathbb{R}^2$. Put $\xi
=(\omega
,\eta
)$ and
%
%e4.31 ###
\begin{eqnarray}\label{gammadef}
\Gamma(s,t):\!&=&\int_s^t\bigl(\xi(u)-\xi(s)\bigr)\otimes d\xi
(u)\nonumber\\
&=&\frac{(\omega(t)-\omega(s))^2}{2}e_1\otimes e_1
+\frac{(\eta(t)-\eta(s))^2}{2}e_2\otimes e_2\nonumber\\[-8pt]\\[-8pt]
&&{}
+\biggl(\int_s^t\bigl(\eta(u)-\eta(s)\bigr)\,d\omega(u)
\biggr)e_2\otimes
e_1\nonumber\\
&&{}
+\biggl(\int_s^t\bigl(\omega(u)-\omega(s)\bigr)\,d\eta(u)
\biggr)e_1\otimes
e_2.\nonumber
\end{eqnarray}
It is easy to check that $\Gamma$ is multiplicative [see the
definition in
\textup{(\ref{multip})}], and we obtain a rough path $(\xi,\Gamma
)$. Moreover,
Theorem \ref{intfbm} enables to consider integrals with respect to
$\xi$,
and, by applying \textup{(\ref{taylor})} and Theorem \ref{intrough},
we see
that they
coincide with the integrals of Appendix \ref{roughpaths}, so they
match the
rough paths theory.

\begin{proposition}
Let $\xi$ be a two-dimensional $H$-fractional Brownian motion for
$1/3<H\le1/2$, and let $\Gamma$ be defined by \textup{(\ref
{gammadef})}. Then the
rough path $(\xi,\Gamma)$ coincides with the rough path constructed by
Coutin and Qian \cite{couqian02} by means of linear interpolation on dyadic
subdivisions.
\end{proposition}

\begin{pf}
It is sufficient to check that the integral $\gamma(s,t)$ of \textup
{(\ref{gammast})}
coincides with the other approach, and actually, we only consider
$\gamma(0,1)=\int_0^1\eta\,d\omega$. For $\omega$ fixed, the integral
$\int\eta\,d\omega$ is in the Gaussian space generated by $\eta$,
so it is
characterized by its covariance with the variables $\eta(t)$. But, for
$\omega$ fixed, $\int\eta\,d\omega^a$ converges in $L^2$ to
$\int\eta\,d\omega$, so
\begin{eqnarray*}
\mathbb{E}^\omega\biggl[\eta(t)\int_0^1\eta\,d\omega\biggr]
&=&\lim_a\int_0^1\mathbb{E}[\eta(t)\eta(s)]\,d\omega
^a(s)\\
&=&\lim_a\int_0^1\bigl(\omega^a(1)-\omega^a(s)\bigr)\frac{\partial}{\partial
s}\mathbb{E}[\eta(t)\eta(s)]\,ds\\
&=&\int_0^1\bigl(\omega(1)-\omega(s)\bigr)\frac{\partial}{\partial
s}\mathbb{E}[\eta(t)\eta(s)]\,ds.
\end{eqnarray*}
Thus the integral is in the closed subspace of $L^2$ generated by the
variables $\omega(u)\eta(t)$, and is characterized by
%
%e4.32 ###
\begin{eqnarray}\label{covar}\qquad\quad
\mathbb{E}\biggl[\omega(u)\eta(t)\int_0^1\eta\,d\omega\biggr]
&=&\int_0^1\mathbb{E}\bigl[\omega(u)\bigl(\omega(1)-\omega(s)\bigr)
\bigr]\frac
{\partial
}{\partial s}
\mathbb{E}[\eta(t)\eta(s)]\,ds\nonumber\\[-8pt]\\[-8pt]
&=&\int_0^1\mathbb{E}[\eta(t)\eta(s)]\frac{\partial
}{\partial s}\mathbb{E}
[\omega(u)\omega(s)]\,ds.\nonumber
\end{eqnarray}
On the other hand, the Coutin--Qian integral $\int\eta\,d_{CQ}\omega$
is also
in this closed subspace, and is characterized by
%
%e4.33 ###
\begin{equation}\label{cqcovar}\qquad
\mathbb{E}\biggl[\omega(u)\eta(t)\int_0^1\eta\,d_{CQ}\omega\biggr]
=\lim_n\int_0^1\mathbb{E}[\eta(t)\eta^n(s)]\frac
{\partial
}{\partial s}
\mathbb{E}[\omega(u)\omega^n(s)]\,ds,
\end{equation}
where $(\omega^n,\eta^n)$ are dyadic approximations of $(\omega,\eta
)$. We
have to prove that the two expressions in \textup{(\ref{covar})} and
\textup{(\ref{cqcovar})}
match. It is clear that the expectations in \textup{(\ref{cqcovar})} converge,
and we
can conclude by standard techniques as soon as we prove that
%
%e4.34 ###
\begin{equation}\label{supsup}
\sup_n\int_0^1\biggl|\frac{\partial}{\partial s}
\mathbb{E}[\omega(u)\omega^n(s)]\biggr|^{1+\varepsilon
}\,ds<\infty
\end{equation}
for some $\varepsilon>0$. But $s\mapsto\mathbb{E}[\omega
(u)\omega
^n(s)]$
is the
dyadic approximation of $s\mapsto\mathbb{E}[\omega(u)\omega
(s)
]$ which
contains two terms ($s^{2H}$ and $|u-s|^{2H}$) depending on $s$ (the term
$u^{2H}$ disappears in the differentiation). If $\{s^{2H}\}^n$ and
$\{|u-s|^{2H}\}^n$ denote their dyadic approximations, then
\[
\biggl|\frac{\partial}{\partial s}\{s^{2H}\}^n\biggr|\le
s^{2H-1},\qquad
\biggl|\frac{\partial}{\partial s}\{|u-s|^{2H}\}^n\biggr|\le|u-s|^{2H-1},
\]
so \textup{(\ref{supsup})} holds provided $(1+\varepsilon)(1-2H)<1$.
\end{pf}

\begin{remark}
It is known from the construction of \cite{couqian02} that the rough path
$(\xi,\Gamma)$ is geometric (it is the limit in $p$-variation of finite
variation paths with their double integrals). However, we do not know
whether it is the limit of $(\omega^a,\eta^a)$ with its double integrals.
\end{remark}

\begin{appendix}

\section*{Appendix}

%s5 ###
\subsection{A mixing property}
\label{mixing}

We give a result about the long-range dependence of increments of a
fractional Brownian motion. This result was used in Proposition~\ref{nala}
but may also be of independent interest. After this work was completed,
a~similar result was proved in \cite{norrossak} with a more functional
analytic method.

\renewcommand{\thetheorem}{A.\arabic{section}}
\begin{theoremm}\label{correl}
Consider a fractional Brownian motion $(W_t;t\in\mathbb{R})$ with parameter
$0<H<1$; for $-\infty\le s\le t\le+\infty$, denote by $\mathcal
{F}_t^s$ the
$\sigma$-algebra generated by the increments $W_v-W_u$, $s<u\le v<t$. Let
$t_0<t_1<t_2$, let $F$ and $G$ be real variables which are respectively
measurable with respect to $\mathcal{F}_{t_0}^{-\infty}$ and
$\mathcal{F}_{t_2}^{t_1}$, and let $q>1$. We suppose that $F$ and $G$
are in
$L^q$. Let
\[
R(t_0,t_1,t_2):=\biggl(\frac{t_2-t_1}{t_1-t_0}\biggr)^{1-H}.
\]
Then if $R(t_0,t_1,t_2)$ is small enough, the product $F G$ is integrable
and
\[
|\mathbb{E}[FG]-\mathbb{E}[F]\mathbb{E}[G]|\le C\|
F\|_q\|
G
\|_q
R(t_0,t_1,t_2)
\]
for some $C=C(q,H)$.
\end{theoremm}

In particular, for $q=2$, we get an upper bound for the correlation
coefficient
\[
\rho(\mathcal{F}_{t_0}^{-\infty},\mathcal{F}_{t_2}^{t_3})
:=\sup\biggl\{\frac{|\operatorname{cov}(F,G)|}{\sqrt
{\operatorname{var}(F) \operatorname{var}
(G)}};
F\in\mathcal{F}_{t_0}^{-\infty}, G\in\mathcal{F}_{t_2}^{t_1}
\biggr\}.
\]
This bound is valid if $R(t_0,t_1,t_2)$ is small enough, but the coefficient
is of course bounded by 1 everywhere. We deduce that the mixing property
%
%e5.1 ###
\setcounter{equation}{0}
\begin{equation}\label{rhofjfk}
\rho\bigl(\mathcal{F}_{j\delta}^{-\infty},\mathcal
{F}_{(k+1)\delta
}^{k\delta
}\bigr)
\le\frac{C_H}{1+|k-j|^{1-H}}
\end{equation}
holds for any $\delta>0$ and any integers $j\le k$ in $\mathbb{Z}$.

\begin{remarkk}
The order of magnitude claimed in the theorem is optimal, as it can be seen
by taking for $F$ and $G$ some increments of $W$. However, in
Proposition~\ref{nala}, we do not use the whole $\sigma$-algebra
$\mathcal{F}_{j\delta}^{-\infty}$, but only $\mathcal{F}_{j\delta
}^{(j-1)\delta}$;
in this case, our estimate is rough but sufficient for our result.
\end{remarkk}

\begin{remarkk}
One can consider the similar problem for the $\sigma$-algebra
generated by
$W_u$, $s\le u\le t$, instead of the increments of $W$. This question is
studied in \cite{berkeshor99}, but the result proved there is not
sufficient for us.
\end{remarkk}

For the proof of Theorem \ref{correl}, let us first introduce some notation
concerning fractional calculus. The fractional integral operator (or
left-sided Riemann--Liouville operator) of order $\alpha>0$ is defined by
\[
I^\alpha
g(t):=\frac{1}{\Gamma(\alpha)}\int_0^t(t-s)^{\alpha-1}g(s)\,ds.
\]
It satisfies $I^{\alpha+\beta}=I^\alpha I^\beta$, and it coincides
with the
iterated integral of $g$ if $\alpha$ is an integer. Moreover $I^\alpha
$ maps
the space $L^q([0,T])$ into itself, and
%
%e5.2 ###
\begin{equation}\label{ialphaphi}
I^\alpha\phi_\beta=\phi_{\alpha+\beta} \qquad\mbox{for
$\phi_\beta(t)=t^\beta/\Gamma(\beta+1)$, $\beta>-1$.}
\end{equation}

Consider the fractional Brownian motion $W$ of the theorem. The result is
trivial if $H=1/2$, so we suppose $H\ne1/2$. From the shift
invariance, we
can also suppose $t_0=0$. The Mandelbrot--Van Ness definition states
that if
$(B_t;t\in\mathbb{R})$ is a double standard Brownian motion, then
%
%e5.3 ###
\begin{equation}\label{repres}
W_t:=C\int_{\mathbb{R}}\bigl((t-s)_+^{H-1/2}-(-s)_+^{H-1/2}\bigr)\,dB_s
\end{equation}
is a fractional Brownian motion for $C>0$; we will choose the normalization
\[
C=C_H:=\Gamma(H+1/2)^{-1}.
\]
We also consider an independent standard Brownian motion $(\overline{B}
_t;t\le0)$,
and we let $\mathcal{F}_0$ and $\mathcal{F}_0'$ be the $\sigma
$-algebras generated
respectively by $(B_s; s\le0)$ and $(B_s,\overline{B}_s; s\le0)$.

\begin{lemmaa}\label{cameron}
Let $(f(t);0\le t\le t_2-t_1)$ be a random function which is
measurable with respect to $\mathcal{F}_0'$ and such that $f(0)=0$. We suppose
that $f=I^{H+1/2}g$ for a function $g$ in $L^2([0,t_2-t_1])$. Consider the
perturbed process
%
%e5.4 ###
\begin{equation}\label{wtilde}
\widetilde{W}_t:=W_t+f(t-t_1) 1_{\{t\ge t_1\}},\qquad  t\le t_2.
\end{equation}
Then, if $G(W)$ is a functional depending (as in Theorem \ref{correl}) on
the increments of $W$ between times $t_1$ and $t_2$,
\begin{eqnarray*}
&&|\mathbb{E}[G(W)\mid\mathcal{F}_0]-
\mathbb{E}[G(\widetilde{W})\mid\mathcal{F}_0]|\\
&&\qquad\le C \mathbb{E}[|G(\widetilde{W})|^q\mid\mathcal{F}_0]^{1/q}
\mathbb{E}[(L^{1/2}e^{C L})^p
\mid\mathcal{F}_0]^{1/p}
\end{eqnarray*}
for $1/p+1/q=1$ and some $C=C(q)$, and with
\[
L:=\int_0^{t_2-t_1}g(s)^2\,ds.
\]
\end{lemmaa}

\begin{pf}
By definition, we have
\[
f(t-t_1) 1_{\{t\ge t_1\}}=C_H\int_{t_1}^{t\vee t_1}(t-s)^{H-1/2}g(s-t_1)\,ds,
\]
so
\[
\widetilde{W}_t=C_H\int_{\mathbb{R}}
\bigl((t-s)_+^{H-1/2}-(-s)_+^{H-1/2}
\bigr)\,d\widetilde{B}_s
\]
with
\[
\widetilde{B}_t=B_t+\int_{t_1}^{t\vee t_1}g(s-t_1)\,ds.
\]
The process $B$ is perturbed after time $t_1$ by an absolutely continuous
process which is $\mathcal{F}_0'$-measurable, so by writing the Cameron--Martin
theorem conditionally on~$\mathcal{F}_0'$,
\begin{eqnarray*}
&&\mathbb{E}[G(W)\mid\mathcal{F}_0']\\
&&\qquad=\mathbb{E}\biggl[G(\widetilde{W})\exp\biggl(-\int_{t_1}^{t_2}g(s-t_1)\,dB_s
-\tfrac12\int_{t_1}^{t_2}g(s-t_1)^2\,ds\biggr)\Bigm|\mathcal
{F}_0'\biggr].
\end{eqnarray*}
By conditioning on $\mathcal{F}_0\subset\mathcal{F}_0'$,
\[
\mathbb{E}[G(W)\mid\mathcal{F}_0]-
\mathbb{E}[G(\widetilde{W})\mid\mathcal{F}_0]
=\mathbb{E}\bigl[G(\widetilde{W})\bigl(\exp(\cdots
)-1\bigr)\bigm|\mathcal{F}
_0\bigr]
\]
and the result follows from H\"{o}lder's inequality and standard estimates on
the moments of $\exp(\cdots)-1$.
\end{pf}

\begin{pf*}{Proof of Theorem \protect\ref{correl}}
We use previous notation, and in particular suppose $H\ne1/2$ and $t_0=0$.
Define
%
%e5.5 ###
\begin{equation}\label{ft}
f(t):=C_H\int_{-\infty}^0
\bigl((t+t_1-s)^{H-1/2}-(t_1-s)^{H-1/2}
\bigr)(d\overline{B}_s-dB_s)
\end{equation}
for $0\le t\le t_2-t_1$. Let us assume that $f$ satisfies the
assumption of
Lemma \ref{cameron} (this will be proved later). Consider the process
$\widetilde{W}$ of \textup{(\ref{wtilde})}, and the process
$\overline{W}$ obtained from
$W$ by
replacing $B$ by $\overline{B}$ on $(-\infty,0]$ in \textup{(\ref
{repres})}, so that
\[
\overline{W}_t=C_H\int_{-\infty}^{t\wedge0}
\bigl((t-s)^{H-1/2}-(-s)^{H-1/2}\bigr)\,d\overline{B}_s
+C_H\int_{t\wedge0}^t(t-s)^{H-1/2}\,dB_s.
\]
Then $\overline{W}$ has the same law as $W$, is independent from
$\mathcal{F}
_0$, and
\[
\widetilde{W}_{t+t_1}-\widetilde
{W}_{t_1}=W_{t+t_1}-W_{t_1}+f(t)=\overline{W}
_{t+t_1}-\overline{W}_{t_1},
\]
so $G(\widetilde{W})=G(\overline{W})$ is independent from $\mathcal
{F}_0$ and has the same
law as $G=G(W)$. Thus we can use
\[
\mathbb{E}[G(\widetilde{W})\mid\mathcal{F}_0]=\mathbb
{E}[G],\qquad
\mathbb{E}[|G(\widetilde{W})|^q\mid\mathcal
{F}_0]^{1/q}=\|
G\|_q
\]
in Lemma \ref{cameron}, so that
\[
|\mathbb{E}[G\mid\mathcal{F}_0]-
\mathbb{E}[G]|
\le C \|G\|_q\mathbb{E}[(L^{1/2}e^{C L})^p
\mid\mathcal{F}_0]^{1/p}.
\]
Thus
\begin{eqnarray*}
|\operatorname{cov}(F,G)|
&\le&
C \|G\|_q \mathbb{E}[|F| \mathbb{E}[
(L^{1/2}e^{C L})^p
\mid\mathcal{F}_0]^{1/p}]\\
&\le& C \|G\|_q \|F\|_q \|L^{1/2}e^{C L}\|_p.
\end{eqnarray*}
In order to conclude, we have to estimate this $L^p$ norm. The formula
\textup{(\ref{ft})} for $f$ can be differentiated, so $f$ is smooth and
\[
f^{(k)}(t)=\frac{1}{\Gamma(H-k+1/2)}
\int_{-\infty}^0(t+t_1-s)^{H-k-1/2}(d\overline{B}_s-dB_s)
\]
for $k\ge1$. In particular,
\[
\|f^{(k)}(t)\|_r=C(t+t_1)^{H-k}
\]
for any $r$ and some $C=C(r,k,H)$. On the other hand [recall the definition
of $\phi_\beta$ in \textup{(\ref{ialphaphi})}],
\[
f=f'(0)\phi_1+I^2(f'')=I^{H+1/2}g
\]
for
\[
g=f'(0)\phi_{1/2-H}+I^{3/2-H}(f'').
\]
In particular, $f$ satisfies the assumption of Lemma \ref{cameron}.
Moreover,
\[
\|f'(0)\|_r \phi_{1/2-H}(t)=C t_1^{H-1}t^{1/2-H},
\]
and
\begin{eqnarray*}
\|I^{3/2-H}(f'')(t)\|_r
&\le& C\int_0^t(t-s)^{1/2-H}\|f''(s)\|_r\,ds\\
&\le& C'\int_0^t(t-s)^{1/2-H}(t_1+s)^{H-2}\,ds\\
&\le& C'' t_1^{H-1}t^{1/2-H}
\end{eqnarray*}
where the last estimate is easily obtained by considering separately the
integrals on $[0,t/2]$ and $[t/2,t]$. Thus we have obtained an estimate for
$\|g(t)\|_r$, and we deduce that
\[
\|L^{1/2}\|_r\le
C \bigl((t_2-t_1)/t_1\bigr)^{1-H}=C R(0,t_1,t_2)
\]
for any $r$ and some $C=C(r,H)$. We still have to prove that the
moments of
$\exp(L)$ are bounded; but, from Jensen's inequality,
\[
\exp(r L)\le\frac1{\mathbb{E}L}
\int_0^{t_2-t_1}\exp\biggl(r \frac{g(s)^2}{\mathbb{E}
g(s)^2}\mathbb{E}L\biggr)\mathbb{E}[g(s)^2]\,ds,
\]
so, since $g(s)$ is Gaussian, this expression has bounded expectation
provided $r \mathbb{E}L<1/2$, and therefore if $R(t_0,t_1,t_2)$ is
small enough.
\end{pf*}

%s6 ###
\subsection{Rough paths}
\label{roughpaths}

Our aim is to describe a part of the rough paths theory through a point of
view which is well adapted to our approach (Theorem \ref{intfbm}). Our
result (Theorem \ref{intrough} below) is in particular comparable to
\cite{lejay03,gubin04,feyelprad06}, and we include for completeness a short
proof which is sufficient for our purpose. Let $\xi(t)$ be a path with
finite $p$-variation, for $p<3$. In this case, we learn from the theory of
rough paths that $\xi$ is not sufficient for the construction of an integral
calculus, but we also need its double integrals. More precisely, let
$\xi(t)$ and $\Gamma(s,t)$ take their values respectively in
$\mathbb{R}
^d$ and
$\mathbb{R}^d\otimes\mathbb{R}^d$. We suppose that
%
%e6.1 ###
\begin{equation}\label{xixi}
|\xi(t)-\xi(s)|\le\mu(t-s)^r,\qquad
|\Gamma(s,t)|\le\mu(t-s)^{2r}
\end{equation}
for $r=1/p$ (continuous paths with finite $p$-variation can be reduced to
this case by a change of time). The path is supposed to be
multiplicative in
the sense
%
%e6.2 ###
\begin{equation}\label{multip}
\Gamma(s,t)=\Gamma(s,u)+\Gamma(u,t)+\bigl(\xi(u)-\xi(s)
\bigr)\otimes
\bigl(\xi(t)-\xi(u)\bigr)
\end{equation}
for $s\le u\le t$. If $r>1/2$, then $\Gamma$ is necessarily the Young
integral
\[
\Gamma(s,t)=\int_s^t\bigl(\xi(u)-\xi(s)\bigr)\otimes d\xi(u),
\]
but if $1/3<r\le1/2$, the function $\Gamma$, when it exists, is not unique;
one can add to it $\phi(t)-\phi(s)$ for any $(2r)$-H\"{o}lder continuous
$\phi$.
Let us now explain how one can define integrals $\int\rho\,d\xi$, in
a way
which coincides with the tree approach of Theorem \ref{intfbm}.

\begin{theoremm}\label{intrough}
Consider paths $(\xi,\Gamma)$ satisfying \textup{(\ref{xixi})} and
$\textup{(\ref{multip})}$,
$\rho$ with values in $\mathcal{L}(\mathbb{R}^d,\mathbb{R}^n)$ (the
space of
linear maps),
and $\Phi$ with values in the space
$\mathcal{L}(\mathbb{R}^d,\mathcal{L}(\mathbb{R}^d,\mathbb
{R}^n))=\mathcal{L}(\mathbb{R}
^d\otimes\mathbb{R}
^d,\mathbb{R}^n)$.
We suppose that
\[
\bigl|\rho(t)-\rho(s)-\Phi(s)\bigl(\xi(t)-\xi(s)\bigr)
\bigr|\le
\mu'(t-s)^{2r}
\]
and
\[
|\Phi(t)-\Phi(s)|\le\mu'(t-s)^r.
\]
For any $s<t$ and any subdivision $\Sigma=(t_k)$ of $[s,t]$, put
%
%e6.3 ###
\begin{equation}\label{gsigma}
g(\Sigma):=\sum_k\bigl(\rho(t_k)\bigl(\xi(t_{k+1})-\xi(t_k)
\bigr)+\Phi(t_k)
\Gamma(t_k,t_{k+1})\bigr).
\end{equation}
Then $g(\Sigma)$ converges as $\max(t_{k+1}-t_k)$ tends to 0, and the limit
$\int_s^t\rho\,d\xi$ satisfies
%
%e6.4 ###
\begin{equation}\label{intrhoxi}
\biggl|\int_s^t\rho\,d\xi-\rho(s)\bigl(\xi(t)-\xi(s)\bigr)
-\Phi(s)\Gamma(s,t)\biggr|\le C \mu \mu' (t-s)^{3r}
\end{equation}
for some $C=C(r)$.
\end{theoremm}

\begin{remarkk}
The identification
$\mathcal{L}(\mathbb{R}^d,\mathcal{L}(\mathbb{R}^d,\mathbb
{R}^n))=\mathcal{L}(\mathbb{R}
^d\otimes\mathbb{R}
^d,\mathbb{R}^n)$
is made through $[G(x)](y)=G(x\otimes y)$.
\end{remarkk}

\begin{remarkk}
We use the simple notation $\int\rho\,d\xi$ though the integral actually
depends on $(\rho,\Phi)$ and $(\xi,\Gamma)$. Notice, however, that if
\[
\limsup_{t\downarrow s}|\xi(t)-\xi(s)|/(t-s)^{2r}=+\infty
\]
for almost any $s$ (and this is the case for an $H$-fractional Brownian
motion and $1/3<r<H\le1/2$), then $\Phi$ is uniquely determined by
$\rho$.
\end{remarkk}

\begin{pf*}{Proof of Theorem \protect\ref{intrough}}
In the proof we will use the following result taken from Young integration.
Let $g(\Sigma)$ be a function defined on finite subdivisions $\Sigma=(t_k)$
of $[s,t]$ and let $\Sigma_k$ be the subdivision with $t_k$ removed. We
suppose that
%
%e6.5 ###
\begin{equation}\label{hypgsigma}
|g(\Sigma)-g(\Sigma_k)|\le C_g(t_{k+1}-t_{k-1})^\kappa
\end{equation}
for some $\kappa>1$. Then $g(\Sigma)$ converges as the mesh of
$\Sigma$
tends to 0, and
\[
|\lim g-g(o)|\le C(\kappa)C_g(t-s)^\kappa
\]
where the trivial subdivision $o=(s,t)$. Let $g$ be the functional of
\textup{(\ref{gsigma})}. Then
\begin{eqnarray*}
g(\Sigma)-g(\Sigma_k)&=&
\rho(t_{k-1})\bigl(\xi(t_k)-\xi(t_{k-1})\bigr)+\Phi
(t_{k-1})\Gamma
(t_{k-1},t_k)\\
&&{} +\rho(t_k)\bigl(\xi(t_{k+1})-\xi(t_k)\bigr)+\Phi
(t_k)\Gamma
(t_k,t_{k+1})\\
&&{} -\rho(t_{k-1})\bigl(\xi(t_{k+1})-\xi(t_{k-1})\bigr)-\Phi
(t_{k-1})\Gamma(t_{k-1},t_{k+1})\\
&=&\bigl(\rho(t_k)-\rho(t_{k-1})\bigr)\bigl(\xi(t_{k+1})-\xi
(t_k)
\bigr)\\
&&{} -\Phi(t_{k-1})\bigl(\xi(t_k)-\xi(t_{k-1})\bigr)\otimes
\bigl(\xi
(t_{k+1})-\xi(t_k)\bigr)\\
&&{} +\bigl(\Phi(t_k)-\Phi(t_{k-1})\bigr)\Gamma(t_k,t_{k+1})
\end{eqnarray*}
where we have used the multiplicative property of $\Gamma$. The condition
\textup{(\ref{hypgsigma})} is satisfied with $\kappa=3r$ and $C_g=2
\mu \mu
'$, so
the result is proved.
\end{pf*}

In particular, we can compute the integral $\int f(\xi)\,d\xi$ of a one-form
by considering $\rho=f(\xi)$ and $\Phi=f'(\xi)$; the property
\textup{(\ref{intrhoxi})} implies that the integral is the limit of
generalized
Riemann sums, so it coincides with the standard rough paths approach.
\end{appendix}

%imsgetref loaded by lrinkeviciute, Friday, June, 16:14:42 2008
%  \lower1.5ex\hbox{`}\hidewidth\crcr\unhbox0}}}

%
\printaddresses


\begin{thebibliography}{37}

%b1 ###
\bibitem{aldous93}
\begin{barticle}[msn]
\bauthor{\bsnm{Aldous},~\bfnm{David}\binits{D.}}
(\byear{1993}).
\btitle{The continuum random tree. {III}}.
\bjournal{Ann. Probab.}
\bvolume{21}
\bpages{248--289}.
%\bmrnumber{
\MR{1207226}
\end{barticle}
\endbibitem

%b2 ###
\bibitem{berkeshor99}
\begin{barticle}[msn]
\bauthor{\bsnm{Berkes},~\bfnm{Istv{\'a}n}\binits{I.}} \AND
  \bauthor{\bsnm{Horv{\'a}th},~\bfnm{Lajos}\binits{L.}}
(\byear{1999}).
\btitle{Limit theorems for logarithmic averages of fractional {B}rownian
  motions}.
\bjournal{J. Theoret. Probab.}
\bvolume{12}
\bpages{985--1009}.
%\bmrnumber{
\MR{1729465}
\end{barticle}
\endbibitem

%b3 ###
\bibitem{bertoin96}
\begin{bbook}[msn]
\bauthor{\bsnm{Bertoin},~\bfnm{Jean}\binits{J.}}
(\byear{1996}).
\btitle{L\'evy Processes}.
\bseries{Cambridge Tracts in Mathematics}
\bvolume{121}.
\bpublisher{Cambridge Univ. Press}, \baddress{Cambridge}.
%\bmrnumber{
\MR{1406564}
\end{bbook}
\endbibitem

%b4 ###
\bibitem{bruneau79}
\begin{bincollection}[msn]
\bauthor{\bsnm{Bruneau},~\bfnm{Michel}\binits{M.}}
(\byear{1979}).
\btitle{Sur la {$p$}-variation d'une surmartingale continue}.
In \bbooktitle{S\'eminaire de Probabilit\'es, XIII (Univ. Strasbourg,
  Strasbourg, 1977/78)}.
\bseries{Lecture Notes in Math.}
\bvolume{721}
\bpages{227--232}.
\bpublisher{Springer}, \baddress{Berlin}.
%\bmrnumber{
\MR{544794}
\end{bincollection}
\endbibitem

%b5 ###
\bibitem{couqian02}
\begin{barticle}[msn]
\bauthor{\bsnm{Coutin},~\bfnm{Laure}\binits{L.}} \AND
  \bauthor{\bsnm{Qian},~\bfnm{Zhongmin}\binits{Z.}}
(\byear{2002}).
\btitle{Stochastic analysis, rough path analysis and fractional {B}rownian
  motions}.
\bjournal{Probab. Theory Related Fields}
\bvolume{122}
\bpages{108--140}.
%\bmrnumber{
\MR{1883719}
\end{barticle}
\endbibitem

%b6 ###
\bibitem{duquesne07}
\begin{bmisc}[vtex]
\bauthor{\bsnm{Duquesne},~\bfnm{T.}\binits{T.}}
(\byear{2006}). The coding of compact real trees by real valued functions. Preprint.
\end{bmisc}
\endbibitem

%b7 ###
\bibitem{duqleg02}
\begin{bmisc}[vtex]
\bauthor{\bsnm{Duquesne},~\bfnm{T.}\binits{T.}} \AND
  \bauthor{\bsnm{Le~Gall},~\bfnm{J.-F.}\binits{J.-F.}}
(\byear{2002}). Random trees, {L}\'evy processes and spatial branching processes.
  \textit{Ast\'erisque}~281.
\end{bmisc}
\endbibitem

%b8 ###
\bibitem{duqleg05}
\begin{barticle}[msn]
\bauthor{\bsnm{Duquesne},~\bfnm{Thomas}\binits{T.}} \AND
  \bauthor{\bsnm{Le~Gall},~\bfnm{Jean-Fran{\c{c}}ois}\binits{J.-F.}}
(\byear{2005}).
\btitle{Probabilistic and fractal aspects of {L}\'evy trees}.
\bjournal{Probab. Theory Related Fields}
\bvolume{131}
\bpages{553--603}.
%\bmrnumber{
\MR{2147221}
\end{barticle}
\endbibitem

%b9 ###
\bibitem{duqleg07}
\begin{barticle}[msn]
\bauthor{\bsnm{Duquesne},~\bfnm{Thomas}\binits{T.}} \AND
  \bauthor{\bsnm{Le~Gall},~\bfnm{Jean-Fran{\c{c}}ois}\binits{J.-F.}}
(\byear{2006}).
\btitle{The {H}ausdorff measure of stable trees}.
\bjournal{ALEA Lat. Am. J. Probab. Math. Stat.}
\bvolume{1}
\bpages{393--415 (electronic)}.
%\bmrnumber{
\MR{2291942}
\end{barticle}
\endbibitem

%b10 ###
\bibitem{evans07}
\begin{bbook}[vtex]
\bauthor{\bsnm{Evans},~\bfnm{Steven~N.}\binits{S.~N.}}
(\byear{2008}).
\btitle{Probability and Real Trees}.
\bseries{Lecture Notes in Math.}
\bvolume{1920}.
\bpublisher{Springer}, \baddress{Berlin}.
\bnote{Lectures from the 35th Summer School on Probability Theory held in
  Saint-Flour, July 6--23, 2005}.
%\bmrnumber{
\MR{2351587}
\end{bbook}
\endbibitem

%b11 ###
\bibitem{evanspitwin06}
\begin{barticle}[msn]
\bauthor{\bsnm{Evans},~\bfnm{Steven~N.}\binits{S.~N.}},
  \bauthor{\bsnm{Pitman},~\bfnm{Jim}\binits{J.}} \AND
  \bauthor{\bsnm{Winter},~\bfnm{Anita}\binits{A.}}
(\byear{2006}).
\btitle{Rayleigh processes, real trees, and root growth with re-grafting}.
\bjournal{Probab. Theory Related Fields}
\bvolume{134}
\bpages{81--126}.
%\bmrnumber{
\MR{2221786}
\end{barticle}
\endbibitem

%b12 ###
\bibitem{falconer90}
\begin{bbook}[msn]
\bauthor{\bsnm{Falconer},~\bfnm{Kenneth}\binits{K.}}
(\byear{1990}).
\btitle{Fractal Geometry}.
\bpublisher{John Wiley \& Sons Ltd.}, \baddress{Chichester}.
\bnote{Mathematical foundations and applications}.
%\bmrnumber{
\MR{1102677}
\end{bbook}
\endbibitem

%b13 ###
\bibitem{feyelprad06}
\begin{barticle}[vtex]
\bauthor{\bsnm{Feyel},~\bfnm{Denis}\binits{D.}} \AND
  \bauthor{\bparticle{de~}\bsnm{La~Pradelle},~\bfnm{Arnaud}\binits{A.}}
(\byear{2006}).
\btitle{Curvilinear integrals along enriched paths}.
\bjournal{Electron. J. Probab.}
\bvolume{11}
\bpages{860--892 (electronic)}.
%\bmrnumber{
\MR{2261056}
\end{barticle}
\endbibitem

%b14 ###
\bibitem{gradnouruva05}
\begin{barticle}[msn]
\bauthor{\bsnm{Gradinaru},~\bfnm{Mihai}\binits{M.}},
  \bauthor{\bsnm{Nourdin},~\bfnm{Ivan}\binits{I.}},
  \bauthor{\bsnm{Russo},~\bfnm{Francesco}\binits{F.}} \AND
  \bauthor{\bsnm{Vallois},~\bfnm{Pierre}\binits{P.}}
(\byear{2005}).
\btitle{{$m$}-order integrals and generalized {I}t\^o's formula: The case of a
  fractional {B}rownian motion with any {H}urst index}.
\bjournal{Ann. Inst. H. Poincar\'e Probab. Statist.}
\bvolume{41}
\bpages{781--806}.
%\bmrnumber{
\MR{2144234}
\end{barticle}
\endbibitem

%b15 ###
\bibitem{gubin04}
\begin{barticle}[msn]
\bauthor{\bsnm{Gubinelli},~\bfnm{M.}\binits{M.}}
(\byear{2004}).
\btitle{Controlling rough paths}.
\bjournal{J. Funct. Anal.}
\bvolume{216}
\bpages{86--140}.
%\bmrnumber{
\MR{2091358}
\end{barticle}
\endbibitem

%b16 ###
\bibitem{istaslang97}
\begin{barticle}[msn]
\bauthor{\bsnm{Istas},~\bfnm{Jacques}\binits{J.}} \AND
  \bauthor{\bsnm{Lang},~\bfnm{Gabriel}\binits{G.}}
(\byear{1997}).
\btitle{Quadratic variations and estimation of the local {H}\"older index of a
  {G}aussian process}.
\bjournal{Ann. Inst. H. Poincar\'e Probab. Statist.}
\bvolume{33}
\bpages{407--436}.
%\bmrnumber{
\MR{1465796}
\end{barticle}
\endbibitem

%b17 ###
\bibitem{kesten86}
\begin{barticle}[msn]
\bauthor{\bsnm{Kesten},~\bfnm{Harry}\binits{H.}}
(\byear{1986}).
\btitle{Subdiffusive behavior of random walk on a random cluster}.
\bjournal{Ann. Inst. H. Poincar\'e Probab. Statist.}
\bvolume{22}
\bpages{425--487}.
%\bmrnumber{
 \MR{871905}
\end{barticle}
\endbibitem

%b18 ###
\bibitem{legall91}
\begin{barticle}[msn]
\bauthor{\bsnm{Le~Gall},~\bfnm{Jean-Fran{\c{c}}ois}\binits{J.-F.}}
(\byear{1991}).
\btitle{Brownian excursions, trees and measure-valued branching processes}.
\bjournal{Ann. Probab.}
\bvolume{19}
\bpages{1399--1439}.
%\bmrnumber{
\MR{1127710}
\end{barticle}
\endbibitem

%b19 ###
\bibitem{leglejan98a}
\begin{barticle}[msn]
\bauthor{\bsnm{Le~Gall},~\bfnm{Jean-Francois}\binits{J.-F.}} \AND
  \bauthor{\bsnm{Le~Jan},~\bfnm{Yves}\binits{Y.}}
(\byear{1998}).
\btitle{Branching processes in {L}\'evy processes: The exploration process}.
\bjournal{Ann. Probab.}
\bvolume{26}
\bpages{213--252}.
%\bmrnumber{
\MR{1617047}
\end{barticle}
\endbibitem

%b20 ###
\bibitem{lejay03}
\begin{bincollection}[msn]
\bauthor{\bsnm{Lejay},~\bfnm{Antoine}\binits{A.}}
(\byear{2003}).
\btitle{An introduction to rough paths}.
In \bbooktitle{S\'eminaire de Probabilit\'es XXXVII}.
\bseries{Lecture Notes in Math.}
\bvolume{1832}
\bpages{1--59}.
\bpublisher{Springer}, \baddress{Berlin}.
%\bmrnumber{
\MR{2053040}
\end{bincollection}
\endbibitem

%b21 ###
\bibitem{lyons98}
\begin{barticle}[msn]
\bauthor{\bsnm{Lyons},~\bfnm{Terry~J.}\binits{T.~J.}}
(\byear{1998}).
\btitle{Differential equations driven by rough signals}.
\bjournal{Rev. Mat. Iberoamericana}
\bvolume{14}
\bpages{215--310}.
%\bmrnumber{
\MR{1654527}
\end{barticle}
\endbibitem

%b22 ###
\bibitem{lyonscl07}
\begin{bbook}[vtex]
\bauthor{\bsnm{Lyons},~\bfnm{Terry~J.}\binits{T.~J.}},
  \bauthor{\bsnm{Caruana},~\bfnm{Michael}\binits{M.}} \AND
  \bauthor{\bsnm{L{\'e}vy},~\bfnm{Thierry}\binits{T.}}
(\byear{2007}).
\btitle{Differential Equations Driven by Rough Paths}.
\bseries{Lecture Notes in Math.}
\bvolume{1908}.
\bpublisher{Springer}, \baddress{Berlin}.
\bnote{Lectures from the 34th Summer School on Probability Theory held in
  Saint-Flour, July 6--24, 2004, With an introduction concerning the Summer
  School by Jean Picard}.
%\bmrnumber{
\MR{2314753}
\end{bbook}
\endbibitem

%b23 ###
\bibitem{lyonsqian02}
\begin{bbook}[msn]
\bauthor{\bsnm{Lyons},~\bfnm{Terry}\binits{T.}} \AND
  \bauthor{\bsnm{Qian},~\bfnm{Zhongmin}\binits{Z.}}
(\byear{2002}).
\btitle{System Control and Rough Paths}.
\bseries{Oxford Mathematical Monographs}.
\bpublisher{Oxford Univ. Press}, \baddress{Oxford}.
\bnote{Oxford Science Publications}.
%\bmrnumber{
\MR{2036784}
\end{bbook}
\endbibitem

%b24 ###
\bibitem{molchan99}
\begin{barticle}[msn]
\bauthor{\bsnm{Molchan},~\bfnm{G.~M.}\binits{G.~M.}}
(\byear{1999}).
\btitle{Maximum of a fractional {B}rownian motion: Probabilities of small
  values}.
\bjournal{Comm. Math. Phys.}
\bvolume{205}
\bpages{97--111}.
%\bmrnumber{
\MR{1706900}
\end{barticle}
\endbibitem

%b25 ###
\bibitem{neveu86}
\begin{barticle}[msn]
\bauthor{\bsnm{Neveu},~\bfnm{J.}\binits{J.}}
(\byear{1986}).
\btitle{Erasing a branching tree}.
\bjournal{Adv. in Appl. Probab.}
\banumber{suppl.}
\bpages{101--108}.
%\bmrnumber{
\MR{868511}
\end{barticle}
\endbibitem

%b26 ###
\bibitem{nevpit89}
\begin{bincollection}[vtex]
\bauthor{\bsnm{Neveu},~\bfnm{J.}\binits{J.}} \AND
  \bauthor{\bsnm{Pitman},~\bfnm{J.}\binits{J.}}
(\byear{1989}).
\btitle{Renewal property of the extrema and tree property of the excursion of a
  one-dimensional {B}rownian motion}.
In \bbooktitle{S\'eminaire de Probabilit\'es XXIII}.
\bseries{Lecture Notes in Math.}
\bvolume{1372}
\bpages{239--247}.
\bpublisher{Springer}, \baddress{Berlin}.
%\bmrnumber{
\MR{1022914}
\end{bincollection}
\endbibitem

%b27 ###
\bibitem{norrossak}
\begin{bmisc}[vtex]
\bauthor{\bsnm{Norros},~\bfnm{I.}\binits{I.}} \AND
\bauthor{\bsnm{Saksman},~\bfnm{E.}\binits{E.}}
(\byear{2007}). Local independence of fractional Brownian motion.
Preprint.
\end{bmisc}
\endbibitem

%b28 ###
\bibitem{nourdin07}
\begin{bmisc}[vtex]
\bauthor{\bsnm{Nourdin},~\bfnm{I.}\binits{I.}}
(\byear{2008}). {A simple theory for the study of SDEs driven by a fractional Brownian
  motion, in dimension one}. In \textit{S\'{e}minaire de Probabilit\'{e}s XLI. Lecture Notes in Math.}
  \textbf{1934}. Springer, Berlin.
\end{bmisc}
\endbibitem

%b29 ###
\bibitem{peres99}
\begin{bincollection}[msn]
\bauthor{\bsnm{Peres},~\bfnm{Yuval}\binits{Y.}}
(\byear{1999}).
\btitle{Probability on trees: An introductory climb}.
In \bbooktitle{Lectures on Probability Theory and Statistics (Saint-Flour,
  1997)}.
\bseries{Lecture Notes in Math.}
\bvolume{1717}
\bpages{193--280}.
\bpublisher{Springer}, \baddress{Berlin}.
%\bmrnumber{
\MR{1746302}
\end{bincollection}
\endbibitem

%b30 ###
\bibitem{picard06}
\begin{barticle}[vtex]
\bauthor{\bsnm{Picard},~\bfnm{Jean}\binits{J.}}
(\byear{2006}).
\btitle{Brownian excursions, stochastic integrals, and representation of
  {W}iener functionals}.
\bjournal{Electron. J. Probab.}
\bvolume{11}
\bpages{199--248 (electronic)}.
%\bmrnumber{
\MR{2217815}
\end{barticle}
\endbibitem

%b31 ###
\bibitem{pisierxu88}
\begin{barticle}[msn]
\bauthor{\bsnm{Pisier},~\bfnm{Gilles}\binits{G.}} \AND
  \bauthor{\bsnm{Xu},~\bfnm{Quan~Hua}\binits{Q.~H.}}
(\byear{1988}).
\btitle{The strong {$p$}-variation of martingales and orthogonal series}.
\bjournal{Probab. Theory Related Fields}
\bvolume{77}
\bpages{497--514}.
%\bmrnumber{
\MR{933985}
\end{barticle}
\endbibitem

%b32 ###
\bibitem{pitman06}
\begin{bbook}[vtex]
\bauthor{\bsnm{Pitman},~\bfnm{J.}\binits{J.}}
(\byear{2006}).
\btitle{Combinatorial Stochastic Processes}.
\bseries{Lecture Notes in Math.}
\bvolume{1875}.
\bpublisher{Springer}, \baddress{Berlin}.
\bnote{Lectures from the 32nd Summer School on Probability Theory held in
  Saint-Flour, July 7--24, 2002, With a foreword by Jean Picard}.
%\bmrnumber{
\MR{2245368}
\end{bbook}
\endbibitem

%b33 ###
\bibitem{rusval93}
\begin{barticle}[msn]
\bauthor{\bsnm{Russo},~\bfnm{Francesco}\binits{F.}} \AND
  \bauthor{\bsnm{Vallois},~\bfnm{Pierre}\binits{P.}}
(\byear{1993}).
\btitle{Forward, backward and symmetric stochastic integration}.
\bjournal{Probab. Theory Related Fields}
\bvolume{97}
\bpages{403--421}.
%\bmrnumber{
\MR{1245252}
\end{barticle}
\endbibitem

%b34 ###
\bibitem{stricker79}
\begin{bincollection}[msn]
\bauthor{\bsnm{Stricker},~\bfnm{C.}\binits{C.}}
(\byear{1979}).
\btitle{Sur la {$p$}-variation des surmartingales}.
In \bbooktitle{S\'eminaire de Probabilit\'es, XIII (Univ. Strasbourg,
  Strasbourg, 1977/78)}.
\bseries{Lecture Notes in Math.}
\bvolume{721}
\bpages{233--237}.
\bpublisher{Springer}, \baddress{Berlin}.
%\bmrnumber{
\MR{544795}
\end{bincollection}
\endbibitem

%b35 ###
\bibitem{williams01}
\begin{barticle}[msn]
\bauthor{\bsnm{Williams},~\bfnm{David R.~E.}\binits{D.~R.~E.}}
(\byear{2001}).
\btitle{Path-wise solutions of stochastic differential equations driven by
  {L}\'evy processes}.
\bjournal{Rev. Mat. Iberoamericana}
\bvolume{17}
\bpages{295--329}.
%\bmrnumber{
\MR{1891200}
\end{barticle}
\endbibitem

%b36 ###
\bibitem{young36}
\begin{barticle}[msn]
\bauthor{\bsnm{Young},~\bfnm{L.~C.}\binits{L.~C.}}
(\byear{1936}).
\btitle{An inequality of the {H}\"older type, connected with {S}tieltjes
  integration}.
\bjournal{Acta Math.}
\bvolume{67}
\bpages{251--282}.
%\bmrnumber{
\MR{1555421}
\end{barticle}
\endbibitem

%b37 ###
\bibitem{zahle98}
\begin{barticle}[msn]
\bauthor{\bsnm{Z{\"a}hle},~\bfnm{M.}\binits{M.}}
(\byear{1998}).
\btitle{Integration with respect to fractal functions and stochastic calculus.
  {I}}.
\bjournal{Probab. Theory Related Fields}
\bvolume{111}
\bpages{333--374}.
%\bmrnumber{
\MR{1640795}
\end{barticle}
\endbibitem

\end{thebibliography}
\end{document}